\documentclass{amsart}
\usepackage{amssymb}
\usepackage{amscd}
\usepackage{amsmath}
\usepackage[dvips]{graphicx}
\usepackage{enumerate}

\newtheorem{theo}{Theorem}
\newtheorem{lema}[theo]{Lemma}
\newtheorem{cor}[theo]{Corollary}
\newtheorem{prop}[theo]{Proposition}
\newtheorem{definition}[theo]{Definition}
\newtheorem{remark}[theo]{Remark}
\newtheorem{notation}[theo]{Notation}

\newcommand{\AAA}{{\mathbb{A}}}
\newcommand{\CC}{{\mathbb{C}}}
\newcommand{\NN}{{\mathbb{N}}}
\newcommand{\RR}{{\mathbb{R}}}
\newcommand{\ZZ}{{\mathbb{Z}}}
\newcommand{\KK}{{\mathbb{K}}}
\newcommand{\calA}{{\mathcal{A}}}
\newcommand{\calC}{{\mathcal{C}}}
\newcommand{\calD}{{\mathcal{D}}}
\newcommand{\calE}{{\mathcal{E}}}
\newcommand{\calF}{{\mathcal{F}}}

\newcommand{\calH}{{\mathcal{H}}}
\newcommand{\calJ}{{\mathcal{J}}}

\newcommand{\calM}{{\mathcal{M}}}
\newcommand{\calN}{{\mathcal{N}}}
\newcommand{\calO}{{\mathcal{O}}}

\newcommand{\calS}{{\mathcal{S}}}
\newcommand{\comp}{{\circ}}
\newcommand{\ocn}{{\calO_{\CC^n,O}}}
\newcommand{\mm}{{\mathbf{m}}}

\newcommand{\codim}{{\mathrm{codim}}}
\newcommand{\mult}{{\mathrm{mult}}}
\newcommand{\wideg}{{\mathrm{wideg}}}

\begin{document}
\title[Finite-determinacy for non-isolated singularities]{Topological finite-determinacy of functions with non-isolated singularities}
\author{Javier Fern\'andez de Bobadilla}
\address{Mathematisch Instituut. Universiteit Utrecht. Postbus 80010. 3508TA Utrecht. The Netherlands.}
\email{bobadilla@math.uu.nl}
\thanks{Supported by the Netherlands Organisation for Scientific Research (NWO). Supported by the Spanish MCyT project BFM2001-1448-C02-01}
\date{}
%\commby{editor}
\begin{abstract}
We introduce the concept of topological finite-determinacy for germs of analytic functions within a fixed ideal 
$I$, which provides a notion of topological finite-determinacy of functions with non-isolated singularities. We prove the following statement 
which generalizes classical results of Thom and Varchenko: let $A$ be the complement in the ideal $I$ of the space 
of germs whose topological type remains unchanged under a deformation within the ideal that only modifies sufficiently 
large order terms of the Taylor expansion; then $A$ has infinite codimension in $I$ in a suitable sense. We also prove the
existence of generic topological types of families of germs of $I$ parametrized by an irreducible analytic set.
\end{abstract}
\maketitle

\section{Introduction}

R.~Thom announced in~\cite{Th} his Stabilisation Theorem stating the following: let $J^r(n,m)$ denote the space of 
$r$-jets of germs at the origin of differentiable mappings from $\RR^n$ to $\RR^m$, and $\pi^s_r:J^s(n,m)\to J^r(n,m)$ the
natural projection mapping. Consider $f\in J^r(n,m)$; there exists a positive integer $s$, depending only on $n$, $m$ and 
$r$, and a proper algebraic subset $\Sigma\subset(\pi^s_r)^{-1}(f)$ such that any two germs $g_1$ and $g_2$ whith 
the same $s$-jet belonging to $(\pi^s_r)^{-1}(f)\setminus\Sigma$ have the same topological type. 

Although R.~Thom gave in~\cite{Th} rather detailed ideas for the proof of his theorem, the first complete proof was given in~\cite{Va1},~\cite{Va2} by
A. Varchenko, and followed a completely different line. More in the line of R. Thom's ideas, E.~Looijenga's thesis contains the result in the function case 
(that is, when $m=1$). Later, A.~du~Plessis (see~\cite{Pl}) gave another proof for arbitrary $m$ 
based on Thom's suggestions, using also his own ideas and ideas from Mather (actually both A.~Varchenko and 
A.~du~Plessis gave slightly stronger statements than R.~Thom's). In this paper we are interested in a generalization of these 
results valid in the realm of (complex or real) analytic non-isolated hypersurface singularities. To see what kind of 
properties are desirable let us state A.~Varchenko's results in the setting of complex analytic 
functions: let $J^r(\CC^n,\CC)_O$ be the 
space of $r$-jets of germs of holomorphic functions at the origin $O$ of $\CC^n$; denote by $\calD$ the group 
of germs of biholomorphisms fixing the origin of $\CC^n$; there is a natural action of $\calD$ in $J^r(\CC^n,\CC)_O$ by composition on the right. 

\begin{theo}[Varchenko~\cite{Va1}]
\label{old} 
Let $T\subset J^r(\CC^n,\CC)_O$ be an irreducible algebraic subset. There exists $s\geq r$, and a proper algebraic subset
$A\subset (\pi^s_r)^{-1}(T)$ such that any two germs $f_1$ and $f_2$ whose $s$-jet is in $(\pi^s_r)^{-1}(T)\setminus A$ 
have the same topological type.

Moreover, for each $r\geq 1$ there exists a partition of $J^r(\CC^n,\CC)_O$ into disjoint constructible subsets 
$U^r_0,...,U^r_{k(r)}$, 
invariant by the action of $\calD$, such that:
\begin{enumerate}
\item If $i>0$, any two germs $f_1$ and $f_2$ whose $r$-jet is in $U^r_i$ has the same topological type.
\item The codimension of $U^r_0$ tends to infinity as $r$ increases. 
\end{enumerate}
\end{theo}

The subsets $U^r_i$ can be constructed so that, if $s>r$ and $i>0$ then $(\pi^s_r)^{-1}(U^r_i)$ coincides 
with one of the subsets $U^s_j$ with $j>0$; this enables to 
decompose $\ocn$ as a union of subsets $\{V_i\}_{i\in\ZZ_{\geq 0}}$ such that, for $i>0$, the $V_i$'s are formed by germs of fixed 
topological type determined by their $r(i)$-jet (for a number $r(i)$  only depending
on $i$), and $V_0$ is infinite codimensional in a suitable sense (and therefore 
easily avoidable by deformation).  

Observe that the topological type of any function $f\in\ocn$ with non-isolated 
singularities is not determined by any $r$-jet of it, no matter how big is $r$:
summation of a generic homogeneous polynomial of degree arbitrarily high transforms 
it into a function with an isolated singularity at the origin, whose sheaf of 
vanishig cycles is concentrated at the origin, unlike the sheaf of vanishing cycles 
of $f$. Therefore all the functions defining non-isolated singularities belong to 
the residual set $V_0$ of the decomposition given above, and, consequently, 
Theorem~\ref{old} is only meaningful for the study of isolated singularities. 

The object of this paper is to prove a replacement of Theorem~\ref{old} which is meaningful for the study of non-isolated 
singularities of (complex or real) analytic functions. Our strategy is to 
work with functions belonging to a fixed ideal $I$ instead of the whole space of analytic germs at the origin (for 
example, if we want to study functions which are singular at a line, we can take 
$I$ to be the square of the ideal defining it). In this paper we prove a generalization of Theorem~\ref{old} valid for any ideal of germs of complex or real
analytic functions. 

Working whithin a fixed ideal has been already 
succesful in the study of non-isolated singulaties: generalized versality and 
analytic finite-determinacy, study of the Milnor fibration... (see for example ~\cite{Bo1},~\cite{Jo},~\cite{Pe2},~\cite{Pe3},~\cite{Ne},~\cite{Si1},~\cite{Si2},~\cite{Za}). 
Many of these papers use a generalized morsification method that consists in 
deforming non-isolated singularities 
within a fixed ideal $I$ to get simpler ones, and then study their properties. Up to now this works only when
$I$ has simple geometric properties: it is the square of a complete intersection ideal defining an isolated singularity, or
the analytic space defined by it is low dimensional. To stablish the Morsification Method in general
it is needed the generalization Theorem~\ref{old} provided in this paper (the generalized Morsification Method will appear in~\cite{Bo2}).  

Furthermore, the study of functions within an ideal 
is also relevant for the study of isolated singularities satisfying a fixed ammount 
of conditions (having some fixed tangencies or multiplicites at infinitely near points...)

Unlike in the case of isolated singularities, in the study of functions with 
non-isolated singularities, the interesting phenomenons are not 
concentrated at the origin, but at a neighbourhood at the origin of the singular 
locus; this makes insufficient in practice the straightforward generalization of 
Theorem~\ref{old}, in which the ring $\ocn$ is replaced by the ideal $I$. We will 
need an stronger formulation in which we consider, instead of jets at the 
origin, jets at points ranging in a fixed neighbourhood at the origin: view the ideal
$I$ as the stalk at the origin of a coherent ideal sheaf $\tilde{I}$ defined in a 
neighbourhood $W$ of the origin. Define 
\begin{equation}
\label{jotaerre}
J^r(W,\tilde{I}):=\coprod_{x\in W}\tilde{I}_x/\mm_x^{r+1}\cap\tilde{I}_x,
\end{equation}
where $\tilde{I}_x$ is the stalk of 
$\tilde{I}$ at $x$, and $\mm_x$ the ideal of analytic functions vanishing at $x$.
For any non-negative integer $r$ we will stratify 
the set $J^r(W,\tilde{I})$ in such a way that each stratum is an analytic variety. The generalization of the
subsets $T$, $A$ and $U^r_i$ in Theorem~\ref{old} will be 
analytic subsets of the strata of $J^r(W,\tilde{I})$ satisfying 
analogous properties. 

As we are working on a neighbourhood $W$ of the origin, the right generalization of the invariance of the $U^r_i$'s by the action of
$\calD$ consists in the property that the subsets generalizing the $U^r_i$'s are invariant by the action of diffeomorphisms between open subsets
of $W$ which preserve the ideal sheaf $\tilde{I}$. Actually, what we will prove is the invariance of such subsets by flows within $W$ preserving
$\tilde{I}$; this can be viewed as an infinitesimal version of the invariance by diffeomorphisms, and turns out to be enough for applications.

The structure of the paper is the following: first we stratify
the spaces $J^r(W,\tilde{I})$ for any $r\leq\infty$, show that the strata are analytic varieties in a natural way when $r<\infty$, and a limit of analytic 
varieties when $r=\infty$. We define the class of closed analytic subsets of the strata of $J^\infty(W,\tilde{I})$, their irreducible components and 
codimension. We also introduce the concept of residual subsets of $J^\infty(W,\tilde{I})$, which, roughly speaking, is a subset of infinite codimension.
Then we state a proposition that generalizes the main proposition of~\cite{Va1} to our setting. It roughly states that given any irreducible analytic subset 
of a stratum of $J^\infty(W,\tilde{I})$ there is a proper closed analytic subset (discriminant) of it such that two germs in the same connected component of 
its complement are topologically finite-determined and topologically equivalent. Then, letting residual subsets enter the picture,  we show that such a 
discriminant is unique provided that it satisfies 
a certain minimality condition.  After we define the concept of {\em flow-invariant} subsets of $J^\infty(W,\tilde{I})$, as a 
replacement of the concept of diffeomorphism-invariant subsets of $\calO_{\CC^n,O}$. Then we can formulate the main result of the paper, which generalizes to 
our setting Theorem~\ref{old}, and also shows that the predicted subsets are minimal and unique in a certain sense. In its proof the invariance of the relevant
subsets by flows preserving  the ideal sheaf $\tilde{I}$ gets involved, in contrast with the proof of the invariance by the action of $\calD$ of 
Theorem~\ref{old}, which is straightforward. We leave for the last section the proof of the proposition stating the existence of discriminants. For 
this we have found A.~Varchenko's ideas rather suitable for our setting; nevertheless, a naive generalization of Varchenko's proof collapses in a 
fundamental way in several places (summation of generic functions of high order is needed, and this takes us out of our ideal sheaf).  

\section{The analytic structure of $J^r(W,\tilde{I})$} 

Let the field $\KK$ be either $\RR$ or $\CC$; denote by $\calE_W$ the sheaf of
analytic functions on an open subset $W$ of $\KK^n$ (when $\KK=\CC$ the sheaf 
$\calE_W$ is the sheaf of holomorphic functions $\calO_W$); for any $x\in W$ 
denote by $\mm_x$ the maximal ideal of the stalk of $\calE_{W,x}$.

When $\KK=\RR$ we will work with an special type of real analytic subsets. Consider $\RR^n$ as the subset of $\CC^n$ consisting of the points with real 
coordinates. Let $W\subset\RR^n$ be an open subset; a $C$-{\em analytic subset} of $W$ is a closed real analytic subset $T\subset W$ such that there exists an 
open neighbourhood $W^*$ of $W$ in $\CC^n$ satisfying $W^*\cap\RR^n=W$ and a closed complex analytic subset $T^*\subset W^*$ such that $T^*\cap\RR^n=T$. A 
{\em Zarisky locally closed} $C$-analytic subset of $W$ is the difference between two $C$-analytic subsets of $W$. We use 
Ch.5~of~\cite{Na} as a general survey reference for $C$-analytic subsets; proofs are due to H.~Cartan, F.~Bruhat and H.~Whitney, and can be found in~\cite{Ca},
and, mostly, in~\cite{BW}. For the convenience of the reader we recall the properties that are convenient for us: 

A real analytic subset $T\subset W$ is $C$-analytic if and only if there exist a coherent ideal sheaf in $\calE_W$ whose zero-set is $T$. 
Any real analytic subset is locally $C$-analytic.
Any (possibly infinite) intersection of $C$-analytic subsets is a $C$-analytic subset. Any locally finite union of $C$-analytic subsets is $C$-analytic. 
The inverse image of a $C$-analytic subset by an analytic mapping is a $C$-analytic subset.
A $C$-analytic subset is {\em $C$-irreducible} if it is not the union of two $C$-analytic subsets different from itself (a $C$-irreducible 
$C$-analytic subset need not be irreducible as a real analytic set). Given a $C$-analytic subset $T\subset W$ there exists a unique irredundant, locally 
finite family of $C$-irreducible $C$-analytic subsets whose union is $T$. There is a notion of dimension of $C$-analytic subsets which satisfies the following
properties: let $T'\subset T$ be $C$-analytic subsets of $W$,  where $T$ is $C$-irreducible, then $\dim(T')<\dim(T)$; if $T$ is a $C$-analytic subset of 
dimension $p$, there is a $C$-analytic subset $T'\subset T$ such that $\dim(T')<p$, and $T\setminus T'$ is an analytic manifold of dimension $p$; also, any 
point of $T$ can be aproximated by points where $T$ is an analytic manifold of dimension $p$.
A {\em complexification} of a real $(C)$-analytic subset $T$ is a complex analytic variety $T^*$ such that $T$ is embedded in $T^*$ as a real analitic variety,
and for each point $t\in T$ there is a complex analytic isomorphism $\varphi:U^*\to Y^*$ from a neighbourhood $U^*$ of $t$ in $T^*$ to a closed complex 
analytic subset $Y^*$ of an open subset of $\CC^n$ such that $Y^*\cap\RR^n=\varphi(T\cap U^*)$. Let $T$ be a $C$-analytic subset of $W$ and $T^*$ a 
complexification of $T$. Then $C\subset T$ is a $(C)$-analytic subset of $W$ if and only if there is an open neighbourhood $U^*$ of $T$ in $T^*$ and a closed 
complex analytic subset of $C^*$ of $U$ such that $C^*\cap T=C$. 

All the properties above are standard in the complex analytic setting whenever they make sense. 

We will adopt the following notational convention: when we work simultaneosly with 
$\KK=\RR,\CC$ and we write ($C$)-{\em analytic}, we mean $C$-analytic when $\KK=\RR$, and complex-analytic when $\KK=\CC$; when we write just {\em analytic},
we mean just real-analytic or complex-analytic depending on whether $\KK$ equals $\RR$ or $\CC$. When we say that a $(C)$-analytic subset is irreducible we 
mean $(C)$ irreducible when $\KK=\RR$, and just irreducible when $\KK=\CC$. 

Let $O$ denote the origin of $\KK^n$. Consider an ideal $I\subset\calE_{O}$; let $\{f_1,...,f_m\}$ be a set generators of it. 
Consider an open neighbourhood $W$ of the origin where each of the generators is defined. Then $\{f_1,...,f_m\}$ generate a coherent ideal 
sheaf $\tilde{I}$ whose stalk $\tilde{I}_O$ is equal to $I$.

For any $V\subset W$ we define
\begin{equation}
\label{jotainfinito}
J^r(V,\tilde{I}):=\coprod_{x\in V}\tilde{I}_x/\mm_x^{r+1}\cap\tilde{I}_x\quad\quad J^\infty(V,\tilde{I}):=\coprod_{x\in V}\tilde{I}_x
\end{equation}
If $0\leq r\leq s\leq \infty$ there are obvious projection mappings 
\begin{equation}
\pi^s_r:J^s(V,\tilde{I})\to J^r(V,\tilde{I}),
\end{equation}
For any $r\leq\infty$ there is another natural projection mapping
\begin{equation}
pr_r:J^r(V,\tilde{I})\to V
\end{equation}
whose fibre $J^r(V,\tilde{I})_x$ over a point $x\in V$ is, if $r<\infty$, the vector
space $\tilde{I}_x/\mm_x^{r+1}\cap\tilde{I}_x$, and, if $r=\infty$, the space 
$\tilde{I}_x$.

For any $x\in W$ define the function $\mu_x:\ZZ_{>0}\to\ZZ$ by the formula
\begin{equation}
\label{musubx}
\mu_x(r):=\dim_{\KK}(\tilde{I}_x/\mm_x^{r+1}\cap\tilde{I}_x). 
\end{equation} 
Consider the Hilbert-Samuel function $H_x$ of the $\calE_{W,x}$-module $\calM_x:=\calE_{W,x}/\tilde{I}_x$. As 
\[H_x(r)=\dim_\KK(\calM_x/\mm_x^{r+1}\calM_x)=\dim_\KK(\frac{\calE_{W,x}/\mm_x^{r+1}}{\tilde{I}_x+\mm_x^{r+1}/\mm_x^{r+1}})\]
 we deduce that $\mu_x(r)=\dim_\KK(\calE_{W,x}/\mm_x^{r+1})-H_x(r)$. In~\cite{BM} it is proved the Zarisky analytic upper-semicontinuity of the function $H_x$; 
therefore the function $\mu$ is Zarisky analytic lower-semicontinous. As any real-analytic subset is locally $C$-analytic, shrinking $W$ we can assume that
the subset of $W$ where the function $\mu$ is smaller or equal than a given function is a closed ($C$)-analytic subset.

\begin{definition}
We define the {\em Hilbert-Samuel stratification of} $W$ {\em with respect to} $\tilde{I}$ to be the minimal partition  of 
$W$ such that $\mu_x=\mu_y$ for any two points $x$ and $y$ in the same stratum. The strata are Zarisky locally closed $(C)$-analytic subsets, and we will call 
them the $\tilde{I}$-\emph{strata} of $W$.
\end{definition}

We will use the following notations: consider an analytic function $f$ on an open subset $U\subset\KK^m$. For any $x\in U$ we denote the germ of $f$ at $x$ by 
$f_x$. For any positive integer $r$ we denote the $r$-jet of $f$ at $x$ by $j^rf_x$. Denote by $J^r(M,\KK^m)$ the manifold of $r$-jets of mappings from an 
analytic manifold $M$ to $\KK^m$; it has a natural structure of vector bundle over $M$; given any subset $X$ of $M$ we denote by $J^r(M,\KK^m)_{|X}$ the 
restriction of the bundle to $X$.

For any positive integer $r$ and any $\tilde{I}$-stratum $X$ of $W$ we endow
\begin{equation} 
\label{candidato}
pr_r:J^r(X,\tilde{I})\to X
\end{equation}
with a natural structure of analytic vector bundle as follows: define the $\calE_W$-epimorphism $\varphi:\calE_{W}^m\to\calE_{W}$ by the formula  
$\varphi(g_{1},...,g_{m}):=\sum_{k=1}^m g_kf_k$. Taking $r$-jets for any positive integer $r$, we obtain a mapping
\begin{equation}
\varphi^r:J^r(W,\KK^m)\to J^r(W,\KK),
\end{equation}
given by the formula $\varphi^r(j^rg_{1,x},...,j^rg_{m,x}):=\sum_{k=1}^m j^r(g_{k}f_k)_x$. We observe that $\varphi^r$ is a homomorphism between trivial 
analytic vector bundles over $W$, whose image is canonically identified with $J^r(W,\tilde{I})$. Then the restriction 
\begin{equation}
\varphi^r_{|X}:J^r(W,\KK^k)_{|X}\to J^r(W,\KK)_{|X}
\end{equation}
is a homomorphism of constant rank between trivial analytic vector bundles; therefore its image $J^r(X,\tilde{I})$ has structure of locally trivial 
analytic vector bundle over $X$ of rank $r(X)=\mu_x(r)$ (for $x\in X$ arbitrary) with projection mapping $pr_r$. Furthermore the inclussion
$J^r(X,\tilde{I})\hookrightarrow J^r(W,\KK)_{|X}$ is a monomorphism of analytic vector bundles.  

We denote by $\partial X$ the closed ($C$)-analytic subset of $W$ given by $\overline{X}\setminus X$.
Clearly $J^r(X,\tilde{I})$ is closed analytic in $J^r(W\setminus\partial X,\KK)$. 
Suppose that $\KK=\RR$; we claim that $J^r(X,\tilde{I})$ is actually closed $C$-analytic in $J^r(W\setminus\partial X,\RR)$. View $\RR^n$ as the set 
of points of $\CC^n$ with real coordinates. Locally around any $x\in W$, each generator $f_i$ is given by a convergent power series; therefore there exists an 
open subset $W^*\subset\CC^n$ such that $W^*\cap\RR^n=W$ and complex analytic functions $f^*_1,...,f^*_m$ defined on $W^*$ extending $f_1,...,f_m$. 
Let $\tilde{I}^*$ be the coherent sheaf generated by them. Given any point $x\in X$ we have
\[\tilde{I}_x^*/\mm_x^{r+1}\cap\tilde{I}_x^*\cong (\tilde{I}_x/\mm_x^{r+1}\cap\tilde{I}_x)\otimes_\RR\CC,\]
where $\mm_x$ denotes respectively in each side the maximal ideal of $\calO_{W^*,x}$ and of $\calE_{W,x}$. It follows easily that there is an 
$\tilde{I}^*$-stratum
$X^*$ of $W^*$ satisfying $X^*\cap W=X$; moreover $\partial X^*\cap W=\partial X$.   
The real analytic manifold $J^r(W\setminus\partial X,\RR)$ is naturally embedded as the real part of $J^r(W^*\setminus\partial X^*,\CC)$; as 
\[J^r(X,\tilde{I})=J^r(X^*,\tilde{I}^*)\cap J^r(W\setminus\partial X,\RR)\]
and $J^r(X^*,\tilde{I}^*)$ is closed analytic in $J^r(W^*\setminus\partial X^*,\CC)$ our claim is proved. 

A subset $C$ of $J^\infty(X,\tilde{I})$ is $r$-{\em determined} if it is of the form $C=(\pi^\infty_r)^{-1}(C')$ for a certain subset $C'$ of 
$J^r(X,\tilde{I})$. The {\em determinacy degree} of a subset $C$ of $J^\infty(X,\tilde{I})$ is the minimal integer $r$ such that $C$ is $r$-determined. 

\begin{definition}
A {\em finitely-determined closed ($C$)-analytic subset of} $J^\infty(X,\tilde{I})$ is a $r$-determined subset for a certain integer $s\geq r$, 
such that $\pi^\infty_r(C)$ is a closed ($C$)-analytic subset in $J^r(X,\tilde{I})$. A {\em finitely-determined locally closed ($C$)-analytic subset} 
is the difference between two finitely-determined closed ($C$)-analytic subsets.   
The {\em irreducible (connected) componets} of a $r$-determined (locally) closed ($C$)-analytic subset $C$ are defined to be the 
inverse images by $\pi^\infty_r$ of the irreducible (connected) componets of $\pi^\infty_r(C)$. 
\end{definition}

We endow $J^\infty(X,\tilde{I})$ with the final topology for the family of projections $\{\pi^\infty_r\}_{r\in\NN}$. Then, a family $\{C_j\}_{j\in J}$ of 
finitely-determined subsets of $J^\infty(X,\tilde{I})$ is locally finite if for any $f\in J^\infty(X,\tilde{I})$ there 
exists a positive integer $r$ and a neighbourhood $U$ of $\pi^\infty_r(f)\in J^r(X,\tilde{I})$ such that $(\pi^\infty_r)^{-1}(U)$ meets only finitely $C_j$'s.
Choosing $r$ high enough we can assume that each of the subsets that $(\pi^\infty_r)^{-1}(U)$ meets are $r$-determined; therefore the union 
$\cup_{j\in J}C_j$ looks locally like a finite-determined subset. This motivates

\begin{definition}
$J^\infty(X,\tilde{I})$.
A {\em closed ($C$)-analytic subset} of $J^\infty(X,\tilde{I})$ is the union of a locally finite family of finitely-determined closed ($C$)-analytic subsets
of  $J^\infty(X,\tilde{I})$. A {\em locally closed ($C$)-analytic subset} is the difference 
between two closed ($C$)-analytic subsets.
The set of {\em irreducible components} of a (locally) closed ($C$)-analytic subset of $J^\infty(X,\tilde{I})$ is defined to be the union of the sets of 
irreducible components of the members of the locally finite family that gives rise to it.
\end{definition}

\begin{definition}
Let $C$ be an $r$-determined irreducible (locally) closed ($C$)-analytic subset of $J^\infty(X,\tilde{I})$. Its {\em codimension} 
$\codim(C,J^\infty(X,\tilde{I}))$ is defined to be the codimension of $\pi^\infty_r(C)$ in $J^r(X,\tilde{I})$.
\end{definition}

The above definition does not depend on $r$ because for any two positive integers $s>r$ the mapping $\pi^s_r:J^s(X,\tilde{I})\to J^r(X,\tilde{I})$ is an
affine bundle, and therefore preserves codimension and irreducibility by inverse image.

Consider a possibly infinite filtration  
\begin{equation}
\label{filtracion}
C_1\supset ...\supset C_i\supset ...
\end{equation}
of closed ($C$)-analytic subsets of $J^\infty(X,\tilde{I})$. We say that 
an irreducible component $C_{i,j}$ of $C_{i}$ is {\em stable} if it is an irreducible component of $C_k$ for any $k\geq i$. The component $C_{i,j}$ is called
{\em strongly unstable} if it does not contain any stable component of $C_k$ for any $k\geq i$. Let $\{C_{i,j}\}_{j\in \NN}$ be the set of strongly unstable
irreducible components of $C_i$ (if there is a finite amount we allow repetition in the indexing).  

The intersection of all the closed subsets of the filtration decomposes naturally as the union 
\begin{equation}
\label{decompo}
\bigcap_{i\in\NN}C_i=Z_1\cup Z_2,
\end{equation}
where $Z_1$ is the union of all the stable irreducible components of the $C_i$'s, and 
\[Z_2:=\bigcup_{\{j_i\}_{i\in\NN}\in \NN^\NN}(\bigcap_{i\in\NN}^\infty C_{i,j_i}).\]

Given any positive integer $c$. Let $\{C'_{j}\}_{j\in\calJ_c}$ be the set of all the strongly unstable components of codimension lower bounded by $c$ 
of any of the $C_i$'s. Let $\calJ'_c$ be the indexes corresponding to the components which are maximal by the inclussion relation among the elements of 
$\{C'_{j}\}_{j\in\calJ_c}$. The following inclusion follows easily from the definition of strongly unstable components and the fact that each $C_i$ is 
a locally finite union of irreducible closed ($C$)-analytic subsets:
\[Z_2\subset\bigcup_{j\in\calJ'_c}C'_j.\]
As the family $\{C'_{j}\}_{j\in\calJ'_c}$ is easily checked to be locally finite, the set $Z_2$ is contained in a closed ($C$)-analytic subset with
all of its irreducible components of codimension lower bounded by $c$. This motivates

\begin{definition}
\label{residual}
A closed subset of $J^\infty(X,\tilde{I})$ is {\em residual} if for any positive integer $c$ it is contained in a closed 
($C$)-analytic subset of $U$ with all its irreducible components of codimension at least $c$. 
\end{definition}

\begin{remark}
Let $Z\subset J^\infty(X,\tilde{I})$ admit a decomposition as a union of a closed $(C)$-analytic subset $Z^{(a)}$ and a residual subset $Z^{(r)}$. The subset
$Z^{(a)}$ is uniquely determined by $Z$, and is called the {\em analytic part} of $Z$. The subset $Z^{(r)}$ is uniquely determined if it is minimal among the 
subsets such that $Z=Z^{(a)}\cup Z^{(r)}$; then it is called the {\em residual part} of $Z$. 

Given a filtration like~(\ref{filtracion}), the intersection of all its terms can be decomposed as the union of a closed $(C)$-analytic subset and a residual 
subset. 
\end{remark}

\section{The topological and finite-determinacy discriminants}

To fix ideas we state what we mean by topological finite-determinacy of functions.

\begin{definition}
Two germs of functions $f:(\KK^n,x)\to(\KK,t)$ and $f':(\KK^n,x')\to(\KK,t')$ are said to be {\em topologically right-left equivalent} (we will say 
{\em topologically equivalent to abreviate}) 
if there are germs of homeomorphisms $\phi:(\KK^n,x')\to (\KK^n,x)$ and $\psi:(\KK,t)\to (\KK,t')$ such that 
$\psi\comp f\comp\phi=g$. A function $f\in\tilde{I}_x$ is called topologically $k$-\emph{determined} with respect to $\tilde{I}$ if 
any other $g\in\tilde{I}_x$ with the same $k$-jet is topologically equivalent to $f$. 
\end{definition}

\begin{prop}
\label{Varchenko}
Let $X$ be any $\tilde{I}$-stratum of $W$; let $T\subset J^\infty(\tilde{I},X)$ be a $r$-determined 
locally closed irreducible ($C$)-analytic subset. There exists an $s\geq r$ and a proper $s$-determined closed ($C$)-analytic subset 
$A$ of $T$ such that any two germs $f,g\in T$ whose $s$-jets $\pi^\infty_s(f)$ and $\pi^\infty_s(g)$ are in the same connected 
component of 
\[\pi^\infty_s(T)\setminus \pi^\infty_s(A)\]
are topologically equivalent. 
\end{prop}

The proof of this proposition will be given in the last section of the paper.\\

\noindent
\textbf{Notation}: let $T$ be a subset of $J^r(X,\tilde{I})$ for a certain $r\leq\infty$; given any $x\in W$ we denote by $T_x$ the fibre 
$(pr_r)^{-1}(x)\cap T$ of the restriction of the mapping $pr_r$ to $T$. If $V$ is a subset of $X$ we denote by $T_{|V}$ the intersection 
$(pr_r)^{-1}(V)\cap T$.\\

Let $T\subset J^\infty(\tilde{I},X)$ be a locally closed ($C$)-analytic subset with irreducible components $\{T_j\}_{j\in J}$. 
For each $j\in J$ let $A_j\subset T_j$ be ($C$)-analytic subset predicted by Proposition~\ref{Varchenko}. By the locally finiteness of the family 
$\{T_j\}_{j\in J}$, the union $A:=\cup_{j\in J} A_j$ is closed ($C$)-analytic in $T$. It is easy to check that  any germ in $T\setminus A$ is topologically 
finite-determined and any two germs in the same path connected component of it are topologically equivalent.  

\begin{prop}
\label{minimal}
Let $T$ be a be a locally closed ($C$)-analytic subset of $J^\infty(\tilde{I},X)$. There exists unique subsets $\Gamma\subset\Delta$ of $T$ not containing any 
irreducible component of it with the following properties
\begin{enumerate}[(i)]
\item We have decompositions $\Delta=\Delta^{(a)}\cup\Delta^{(r)}$ and $\Gamma=\Gamma^{(a)}\cup\Gamma^{(r)}$
where $\Delta^{(a)}$ and $\Gamma^{(a)}$ are closed $(C)$-analytic subsets, and 
$\Delta^{(r)}$ and $\Gamma^{(r)}$ are residual closed subsets.
\item Any $f,g\in T$ in the same path-connected component of $T\setminus\Gamma$ are topologically equivalent. Any $f\in T\setminus\Delta$ is moreover 
topologically finite-determined with respect to $\tilde{I}$.
\item The subsets $\Delta$ and $\Gamma$ are minimal among the subsets of $T$ satisfying Properties~(i)~and~(ii).
\end{enumerate}

For any positive integer $k$ we let $T_{\leq k}$ be the union of the irreducible components of $T$ which are $K$-determined, and $T_{>k}$ the union of all the
other irreducible componnets. There is a unique subset $\Delta_k$ of $T_{\leq k}$, not containing any irreducible component of it, with property of beeing 
minimal among the closed subsets $Z$ of $T_{\leq k}$ decomposable as the union of a closed  $(C)$-analytic subset and a residual subset and such that
such that any two germs in the same path-connected component of $T_{\geq k}\setminus Z$ are topologically $k$-determined and have the same topological type. 
The subset $\Delta_k$ is $k$-determined closed $(C)$-analytic. 

We have $\Gamma\subset\Delta\subset\cap_{k\geq s}(\Delta_k\cup T_{>k})$; moreover 
\begin{equation}
\label{igualdad}
\Gamma^{(a)}=\Delta^{(a)}=(\cap_{k\geq s}(\Delta_k\cup T_{>k}))^{(a)},
\end{equation}
in other words, the subsets $\Gamma$, $\Delta$, and the 
intersection $\cap_{k\geq s}(\Delta_k\cup T_{>k})$ only may differ in a residual set.
\end{prop}
\begin{proof}
Let $\calS$ be the set whose elements are subsets of $T$ not containing any of its irreducible components and satisfying the first two properties of $\Gamma$. 
The set $\calS$ is not empty because 
the subset $A$ constructed in the last paragraph before the proposition belongs to it. We consider in it the partial order given by 
inclussion. Consider a chain 
\begin{equation}
\label{chain1}
K_1\supset ...\supset K_i\supset...
\end{equation}
of subsets of $\calS$. We claim that the intersection $K:=\cap_{i\in\NN}K_i$ belongs to $\calC$. 

By Property~(i) each $K_i$ decomposes as 
$K_i:=K_i^{(a)}\cup K_i^{(r)}$, with $K_i^{(a)}$ a closed $(C)$-analytic subset and $K_i^{(r)}$ residual. 
The closed $(C)$-analytic parts form a chain 
\begin{equation}
\label{chain2}
K_1^{(a)}\supset ...\supset K_i^{(a)}\supset...
\end{equation}
We construct a chain 
\begin{equation}
\label{chain3}
L_1\supset ...\supset L_i\supset...
\end{equation}
of subsets that admit a decomposition in a closed $(C)$-analytic subset $L^{(a)}_i$ and a residual subset $L^{(r)}_i$ 
such that $\cap_{i\in\NN}L_i=K$ and all the irreducible components of $L^{(a)}_i$ are either stable or of codimension at least $i$. We proceed 
inductively: suppose that for a certain positive integer $m$ we have defined a chain
\begin{equation}
\label{chain4} 
L_1\supset ...\supset L_m\supset L_{m,m+1}\supset ...\supset L_{m,m+i}\supset ...
\end{equation}
such that, for any $i\leq m$, all the irreducible components of $L^{(a)}_i$ are either stable or of codimension at least $i$ and 
\begin{equation}
\label{intertrunc}
(\bigcap_{i\leq m}L_i)\bigcap(\bigcap_{i>m}L_{m,i})=K.
\end{equation}
Clearly for any $j\in\NN$ there are no non-stable components of codimension strictly smaller than $m$ in $L_{m,m+j}^{(a)}$.
Let $\{C_h\}_{h\in H}$ be the collection of irreducible components of $L_{m,m+1}^{(a)}$ of codimension $m$ which are non-stable in the filtration
given by the closed $(C)$-analytic parts of the elements of filtration~(\ref{chain4}). For any $h\in H$ 
there exists a smallest positive integer $i_h$ such that $C_h$ is not an irreducible component of $L_{m,m+i_h}^{(a)}$. 
For any non-negative integer $j$ let $\{Z_l\}_{l\in L}$ be the set of irreducible components of $L_{m,m+j}$ not belonging to $\{C_h\}_{h\in H}$. Define
\begin{equation}
L_{m+1,m+j}:=(\bigcup_{l\in L}Z_l)\bigcup L_{m,m+j}^{(r)} \bigcup (\bigcup_{j<i_h} (C_h\cap L_{m,m+i_h})).
\end{equation}
Define $L_{m+1}=L_{m+1,m+1}$. By construction, the Equality~(\ref{intertrunc}) holds repalcing $m$ by $m+1$. If we prove 
that each subset $L_{m+1,m+j}$ for $j\geq 0$ admits a decomposition in a closed $(C)$-analytic subset $L^{(a)}_i$ and a residual subset $L^{(r)}_i$, then, by 
construction, all the non-stable irreducible components of $L_{m+1}^{(a)}$ are of codimension at least $m+1$. Iterating the procedure inductively we obtain 
the desired chain~(\ref{chain3}). 

For any non-negative integer $j$ we consider the decomposition 
$L_{m+1,m+j}=L_{m+1,m+j}^{(a)}\cup L_{m+1,m+j}^{(r)}$, where 
\begin{itemize}
\item the set $L_{m+1,m+j}^{(a)}$ is the union of all the irreducible components of $L_{m,m+j}^{(a)}$ not 
belonging to $\{C_h\}_{h\in H}$, together with $\cup_{j<i_h}(C_h\cap L_{m,m+i_h}^{(a)})$. 
\item the set $L_{m+1,m+j}^{(r)}$ is the union of $L_{m,m+j}^{(r)}$ together with $\cup_{j<i_h}(C_h\cap L_{m,m+i_h}^{(r)})$.
\end{itemize}
As the irreducible components $\{C_h\}_{h\in H}$ form a locally finite family, the set $L_{m+1,m+j}^{(a)}$ is a locally finite union of closed $(C)$-analytic 
subsets, and hence it is closed $(C)$-analytic. 
On the other hand, for any positive integer $c$ and any $h$ such that $j<i_h$ 
there exists a closed $(C)$-analytic subset $C'_h$ contained in $C_h$ and containing $C_h\cap L_{m,m+i_h}^{(r)}$ with all its irreducible components of 
codimension at least $c$. By the locally finiteness of $\{C_h\}_{h\in H}$ the subset $\cup_{j<i_h}C'_h$ is closed $(C)$-analytic. 
Therefore $\cup_{j<i_h}(C_h\cap L_{m,m+i_h})^{(r)}$ is residual, and hence $L_{m+1,m+j}^{(r)}$ is also residual.

Now we prove that $K$ belongs to $\calC$. Let $K^{(a)}$ be the union of the stable irreducible components of the filtration given by the $L_i^{(a)}$'s. Define 
$D$ to be the union of all the intersections of the form $\cap_{i\in \NN}(L_i^{(u)}\cup L^{(r)})$ where $L_i^{(u)}$ is the union of all the non-stable 
components of $L_i$. As $K=\cap_{i\in\NN}L_i$ we obtain 
easily the decomposition $K=K^{(a)}\cup D$. We show that $D$ is residual: let $c$ be any positive integer. Consider a closed $(C)$-analytic subset $C$ 
containing $L_c^{(r)}$ with all its irreducible components of codimension at least $c$. The set $L_c^{(u)}\cup C$, whose irreducible components are all of them
of codimension at least $c$, contains $D$. 

Let $\gamma:[0,1]\to T\setminus K$ be a continous path. To show that $K$ belongs to $\calC$ it only remains to chech that the topological type of the germ 
$\gamma(t)$ is independent of $t$. The set $T\setminus K$ is the union of the increasing sequence of open subsets $\{T\setminus K_i\}_{i\in\NN}$. By the 
compactness of $[0,1]$ there is an index so that $\gamma([0,1])\subset T\setminus K_i$. As $K_i\in\calC$ the topological type remains constant along $[0,1]$.
 
We have shown that any decreasing sequence in $\calC$ has a lower bound. By Zorn's Lemma we deduce the existence of $\Gamma$. The uniqueness holds as 
the intersection of two subsets in $\calC$ is easily shown to belong to $\calC$. The existence and uniqueness of $\Delta$ and $\Delta_k$ for any $k\in\NN$ 
is analogous.

Now we show that $\Delta_k$ is $k$-determined $(C)$-analytic. Consider an irreducible component $Z$ of $\Delta^{(a)}$. Let $r$ be the determinacy degree  
of $Z$. If $r\leq k$ then $Z$ is $k$-determined; we study the case $r\geq k$. Consider the affine bundle $\pi^r_k:J^r(X,\tilde{I})\to J^k(X,\tilde{I})$; let 
$B:=\pi^\infty_k(T_{\leq k})$, $E:=\pi^\infty_r(T_{\leq k})$, and $\pi:=\pi^r_{k|E}:E\to B$. The irreducible closed $(C)$-analytic subset 
$Z_r:=\pi^\infty_r(Z)$ is contained in $E$. We claim that the set 
\[Y:=\{y\in B:E_{y}\subset Z_r\}\]
is a (possibly empty) closed $(C)$-analytic subset of $B$. If $\KK=\RR$ then it is easy to choose complexifications $E^*$ and $B^*$ of $E$ and $B$, and a 
mapping $\pi^*:E^*\to B^*$, which is a complex affine bundle, such that $\pi^*_{|E}=\pi$. As $Z_r$ is $C$-analytic there exists an open neigbourhood $U^*$ of 
$E$ in $E^*$ and an irreducible complex closed analytic subset $Z_r^*$ of $U^*$ such that $Z_r^*\cap E=Z_r$. The open subset $U^*$ can be choosen so that for 
any $x\in B$ the fibre $U^*_x$ is connected (we prove this at the end). Let $N$ be the rank of the affine bundle $\pi^*:E^*\to B^*$; the subset
\[A^*:=\{z\in Z^*:\dim_z(Z^*_{\pi^*(z)})=N\}\] 
is a complex closed analytic subset of $Z^*$. Given any $z\in A^*$ we have $\dim_zZ^*_{\pi^*(z)}=\dim_z(E^*_{\pi^*(z)})$. Therefore $Z^*_{\pi^*(z)}$ contains 
an open neighborhood $V_z$ of $z$ in $E^*_{\pi^*(z)}$. Clearly $V_z$ is contained in $A^*_{\pi^*(z)}$, and hence this subset is both closed an open in the 
connected set $U^*\cap E^*_{\pi^*(z)}$. Thus $A^*_{\pi^*(z)}=U^*\cap E^*_{\pi^*(z)}$ for any $z\in A^*$. This implies 
\[U^*\cap (\pi^*)^{-1}(\pi^*(A^*))=A^*,\]
and from here it is easy to deduce that $\pi^*(A^*)$ is a closed complex analytic subset of the open subset $\pi^*(U^*)\subset B^*$. For any $x\in B$ we have 
that $U^*\cap E^*_x\supset E_x$; therefore $\pi^*(A^*)\cap B\subset Y$. On the other hand, if $x\in Y$ we have $E_x\subset U^*\cap Z^*_x$; hence $Z^*_x$ is a 
closed complex analytic subset containing the real part $E_x$ of $U^*\cap E^*_x$; this implies that $Z^*_x$ contains a neighbourhood of $E_x$ in 
$U^*\cap E^*_x$, and hence, by connectedness of $U^*\cap E^*_x$, it is equal to it. Thus $\pi^*(A^*)\cap B=Y$, which proves our claim when $\KK=\RR$. The proof
in the complex case is analogous, but easier. 

If $(\pi^\infty_k)^{-1}(Y)=Z$ then $Z$ is actually $k$-determined. Otherwise we let $\Delta'_k$ be the $k$-determined subset given by the union of 
$\Delta^{(r)}$, $(\pi^\infty_k)^{-1}(Y)$, and all the irreducible components of $\Delta_k^{(a)}$ different from $Z$. Clearly $\Delta'_k$ is strictly contained 
in $\Delta_k$. Consider two germs $f,g$ in the same path-connected component of $T_{\leq k}\setminus\Delta'_k$. We claim that both of them are topologically 
$k$-determined and have the same topological type. This clearly gives a contradiction with the minimality of $\Delta_k$, which shows that $Z$ is 
$k$-determined. In the special case that  neither $f$ nor $g$ does belong to $\Delta_k$ the claim holds by definition of 
$\Delta_k$. Suppose that $f$ belongs to $\Delta_k$. As $f$ does not belong to $\Delta'_k$ there exists an open neighbourhood $V_f$ 
of $f$ in $J^\infty(X,\tilde{I})$ such that $V_f\cap\Delta'_k=\emptyset$. Then $f$ must be an element of $Z\setminus(\pi^\infty_k)^{-1}(Y)$. 
Let $x:=pr_\infty(f)$; as $f$ does not belong to $(\pi^\infty_k)^{-1}(Y)$ there exists a continous path 
$\gamma:[0,1]\to V_f\cap T_{\leq k}\cap J^\infty(X,\tilde{I})_x$ such that $\gamma(0)=f$, the $k$-jet of $\gamma(t)$ does not depend on $t$, and $\gamma(t)$ 
does not belong to $Z$ for $t\neq 0$. Obviously $\gamma(t)\not\in\Delta_k$ for $t\neq 0$, and threfore the claim is reduced to the already verified case in 
which neither $f$ nor $g$ does belong to $\Delta_k$. 

To show the $k$-determinacy of $\Delta_k$ it remains to prove $\Delta_k^{(r)}=\emptyset$, but this follows by arguments analogous to the last paragraph. 
 
The inclussions $\Gamma\subset\Delta\subset\cap_{k\geq s}(\Delta_k\cap T_{>k})$ are trivial. Now we show Equality~(\ref{igualdad}). Let $C$ be a irreducible 
component of $(\cap_{k\geq s}(\Delta_k\cup T_{>k}))^{(a)}$. We suppose that $C$ is not contained in $\Gamma$, and look for a contradiction. There 
exists a certain integer $r$ such that $C$ is an irreducible component of $\Delta_k$ for any $k\geq r$ (this is because 
choosing $r$ large enough we can assume that $T_{>r}$ does not contain $C$); in particular it is $r$-determined.  
By Proposition~\ref{Varchenko} we find an integer $s\geq r$ and a $s$-determined proper closed $(C)$-analytic subset $A_1\subset C$ such 
that any germ in $C\setminus A_1$ is $s$-determined. Define $\Delta'_s$ as the union of all the irreducible components of $\Delta_s$ different from $C$. Let
$A_2$ be the union of all the irreducible components of $\Gamma^{(a)}$ not contained in $\Delta'_s\cup T_{>s}$; as $\Gamma\subset \Delta_s\cup T_{>s}$ we have
$A_2\subset C$, and the inclussion is strict for beeing $C$ not contained in $\Gamma$. Consider a proper closed $(C)$-analytic subset $A_3$ of $C$ such that 
$\Gamma^{(r)}$ is contained in $A_3\cup\Delta'_s\cup T_{>s}$ (the existence of $A_3$ is clear as $\Gamma^{(r)}$ is residual). Define 
\[A:=A_1\cup A_2\cup A_3\cup (C\cap T_{>s}),\]
Any germ in $T_{\leq s}\setminus (A\cup \Delta'_s)$ belongs either to $T_{\leq s}\setminus\Delta_s$ or to $C\setminus A_1$; therefore it is topologically 
$s$-determined. As $T_{\leq s}\setminus (A\cup \Delta'_s)$ is clearly included $T\setminus\Gamma$ any two germs in the same path-connected component of it
have the same topological type. Then $A\cup \Delta'_s$ have the same properties of $\Delta_s$ and is strictly smaller than it. This is a contradiction.

Finally let's check that $U^*$ can be choosen $U^*\cap E^*_x$ is connected for any $x\in B$. As $\pi:E\to B$ has contractible fibres there is a continous
section $s:B\to E$. Therefore, we can give a continous $\RR$-vector bundle structure to $\pi:E\to B$ and $\pi^*:E^*\to B^*$ in such a way that $E$ is a 
subbundle of the restriction of $E^*$ to $B$. 
Shrinking enough $B^*$ around $B$ we can suppose that $B$ is a strong deformation retract of $B^*$ (by the Triangulation Theorem for real analytic subsets 
we can think of $(B^*,B)$ as a polyhedral pair; then we apply Corollary~11~of~\cite{Sp},~page~124). By well known areguments it follows  that there is a vector
subbundle $q:F\to B^*$ of $\pi:E^*\to B^*$ extending $\pi:E\to B$. Using partitions of unity we can construct a continous tensor on $E^*_{|B}$ that restricts 
to an euclidean inner product on each fibre $E_x$; let $d_x$ be the distance induced by it in $E_x$. It is easy to find positive continous 
functions 
\[\alpha:B^*\to\RR\cup\{+\infty\}\quad\quad\alpha:B^*\to\RR,\]
(where a basis of neighbourhoods for $+\infty$ in $\RR\cup\{+\infty\}$ is given by $\{(a,+\infty)\}_{a\in\RR}$), such that $\alpha(B)=\{+\infty\}$ and 
\[\{z\in E^*_x:d_x(z,0)<\alpha(x),\quad d(z,F_x)<\beta(x)\}\subset U^*\cap E^*_x.\]
Redefining  $U^*$ as the open subset 
\[\{z\in E^*:d_{\pi^*(z)}(z,0)<\alpha(\pi^*(z)),\quad d(z,F_{\pi^*(z)})<\beta(\pi^*(z))\}\]
we obtain the desired properties. 
\end{proof}

\begin{definition}
\label{discrimi}
Let $T$ be a be a locally closed ($C$)-analytic subset of $J^\infty(\tilde{I},X)$. 
We call the sets $\Gamma$ and $\Delta$ constructed in the last proposition {\em topological} and {\em finite-determinacy discriminants} of $T$ respectively.
For any positive integer $k$ the set $\Delta_k$ is called the  $k$-\emph{determinacy discriminant} of $T$.
\end{definition}

Next we study the behaviour of discriminants when restricting to open subsets of $U$:

\begin{lema}
\label{restrabierto}
Let $T$, $\Gamma$ and $\Delta$ be as in Lemma~\ref{minimal}. Let $U$ be a open subset of $W$.
When $\KK=\CC$ the subsets $\Gamma_{|U}$ and $T\cap \Delta_{|U}$ are respectively the topological and 
finite-determinacy discriminants of $T_{|U}$. 
\end{lema}
\begin{proof}
Clearly $\Gamma_{|U}$ contains the topological discriminant of $T_{|U}$.
Consider the decomposition $\Gamma_{|U}=\Gamma^{(a)}_{|U}\cup \Gamma^{(r)}_{|U}$, where the first component is $(C)$-analytic and the 
second residual. Any irreducible component $C'$ of $\Gamma^{(a)}_{|U}$ is a subset of a unique irreducible component $C$ of 
$\Gamma^{(a)}$, which is $r$-determined for a certain $r$. By Proposition~\ref{Varchenko} there exists $s\geq r$ and a $s$-determined closed $(C)$-analytic 
subset $A$ of $C$ such that any two 
germs in the same path-connected component of $C\setminus A$ have the same topological type. Let $Z$ be any irreducible component of $T$, as $\KK=\CC$ we have 
that both $Z\setminus\Gamma\cap Z$ and $C\setminus A$ are path-connected, and hence two germs contained in the same of these two subsets have the same 
topological type. We claim that there exists a component $Z$ of $T$ containing $C$ such that the topological type
of the germs of $Z\setminus\Gamma\cap Z$ is different to the topological type of the germs of $C\setminus A$; otherwise we let $\Gamma_1'$ be the union of 
all the irreducible components of $\Gamma^{(a)}$ different from $C$, let $\Gamma'_2$ be the union of the intersections $C\cap X$ where $X$ is any 
irreducible component of $T$ not containing $C$, and define 
\[\Gamma':=\Gamma'_1\cup\Gamma'_2\cup \Gamma^{(r)}\cup A.\]
It is easy to check that $\Gamma'$ has the first two properties of $\Gamma$ and is strictly smaller than it. This is a contradiction.

Let $Z'$ be an irreducible component of $Z_{|U}$ containing $C'$. As 
\[\dim(Z')=\dim(Z)>\dim(C)=\dim(C'),\] 
the inclussion $C'\subset Z'$ is strict.  
As $Z'\setminus \Gamma_{|U}\cap Z'$ and $C'\setminus (A\cap C')$ are contained in $Z\setminus\Gamma\cap Z$ and $C\setminus A$ 
respectively, the topological type of the germs of $Z'\setminus \Gamma_{|U}\cap Z'$ is different from the topological type of the germs of 
$C'\setminus (A\cap C')$. It follows that $C'$ is contained in the topological discriminant of $T_{|U}$.

The proof for $\Delta$ is analogous.
\end{proof}

\section{The main result} 

For any $x\in W$ we denote by $\calD_x$ the group of germs of analytic diffeomorphisms of $\KK^n$ fixing $x$, and by 
$\calD_{x,e}$ the set of germs of analytic diffeomorphisms at $x$ that not necessarily fix it. Following~\cite{Pe2} we 
define $\calD_{\tilde{I}_x,e}$ to be the subset of $\calD_{x,e}$ preserving the ideal; that is, the subset formed by the germs having a representative
$\phi:U\to W$ such that $\phi^*\tilde{I}_{\phi(y)}=\tilde{I}_y$ for any $y\in U$. 
          
Let $\phi_t$ be a 1-parameter family of diffeomorphisms of $\calD_{\tilde{I}_x,e}$ smoothly depending on $t$, such 
that $\phi_0=\mathrm{Id}_{(\KK^n,x)}$; let $\phi_{1,t},...,\phi_{n,t}$ be its components. The germ (at $x$) of analytic vector 
field defined by $X:=\sum_{i=1}^td\phi_{i,t}/dt_{|t=0}\partial/\partial x_i$ preserves the ideal sheaf $\tilde{I}$, that is, satisfies 
$X(\tilde{I}_x)\subset \tilde{I}_x$. 
Let $\Theta$ be the sheaf of analytic vector fields
in $W$; define $\Theta_{\tilde{I},e}$ as the subsheaf whose sections preserve the ideal sheaf $\tilde{I}$; denote the stalk of 
$\Theta_{\tilde{I},e}$ at $x$ by $\Theta_{\tilde{I}_x,e}$. 
Integration associates to any $X\in\Theta_x$, a 1-parameter flow $\phi_t$ of germs of analytic diffeomorphisms of $\calD_{x,e}$ for which 
$\phi_0=\mathrm{Id}_{(\KK^n,x)}$, and such that, if $X\in\Theta_{\tilde{I}_x,e}$, then $\phi_t\in\calD_{\tilde{I}_x,e}$ for any value of $t$. 

Any representative $\phi:U\to W$ of a germ $\phi\in\calD_{I_x,e}$ induces by pushforward a mapping 
\begin{equation}
\phi_*:J^\infty(U,\tilde{I})\to J^\infty(\phi(U),\tilde{I}), 
\end{equation}
defined by $\phi_*(f_y):=(f_y\comp\phi^{-1})_{\phi(y)}$ for any $y\in U$ and $f_y\in\tilde{I}_y$.
As the $\mm_{y}$-adic filtration is transformed by pushforward into the $\mm_{\phi(y)}$-adic filtration the mapping $\phi_*$ descends to a 
mapping
\begin{equation}
j^r\phi_*:J^r(U,\tilde{I})\to J^r(\phi(U),\tilde{I}).
\end{equation}
Clearly any representative $\phi:U\to W$ of a germ $\phi\in\calD_{I_x,e}$ preserves the Hilbert-Samuel stratification (that is $\phi(U\cap X)=\phi(U)\cap X$ 
for any $\tilde{I}$-stratum $X$). It is easy to check that the restriction 
\begin{equation}
\label{analitica}
j^r\phi_*:J^r(U\cap X,\tilde{I})\to J^r(\phi(U)\cap X,\tilde{I})
\end{equation}
is an analytic diffeomorphism when $r<\infty$.

\begin{definition}
\label{flowinv}
Let $X$ be a $\tilde{I}$-stratum of $W$ , and $T\subset J^\infty(X,\tilde{I})$ be a (locally) closed ($C$)-analytic subset. We say that $T$ is 
{\em flow-invariant} if for any open subset $V\subset W$, any vector field $\theta\in\Gamma(V,\Theta_{\tilde{I},e})$ and any flow $\phi:U\times (a,b)\to V$ 
integrating $\theta$, given any $t\in (a,b)$ and $x\in U$, we have $\phi_{t*}(T_x)=T_{\phi_t(x)}$.
\end{definition}

Now we are ready to state the main result of the paper:\\

\noindent
\textbf{Main Theorem.} {\em
Shrink $W$ so that $\Theta_{\tilde{I},e}$ is generated by sections defined on the whole $W$. Consider $X$, a $\tilde{I}$-statum of $W$. Let $T$ be any locally 
closed ($C$)-analytic subset of $J^{\infty}(X,\tilde{I})$ there exist a unique filtration (which we call the {\em filtration by succesive discriminants}) 
\[T=A_0\supset A_1\supset ...\supset A_i\supset ...\]
by closed $(C)$-analytic subsets, and two residual subsets $\Gamma^{(r)}$ and $\Delta^{(r)}$ (called respectively the {\em topological and 
finite-determinacy cummulative residual discriminants} of $T$),  with the 
following properties: 
\begin{enumerate}
\item We have $\cap_{i\geq 0}A_i\subset\Gamma^{(r)}\subset\Delta^{(r)}$. 
\item For any $i\geq 0$ the sets $A_{i+1}\cup (\Gamma^{(r)}\cap A_i)$ and $A_{i+1}\cup (\Delta^{(r)}\cap A_i)$ are 
respectively the topological and finite-determinacy discriminants of $A_i$.
\item Any irreducible component of $A_i$ has codimension at least $i$.  
\item If $T$ is flow-invariant then $A_i$ is flow-invariant for any $i\geq 0$. Moreover the $k$-determinacy
discriminant of $A_i$ is flow invariant. Therefore the sets $\Gamma^{(r)}\cap A_i$ and $\Gamma^{(r)}\cap A_i$ are contained in a 
residual subset which is an intersection of flow-invariant closed $(C)$-analytic subsets of $C$.
\end{enumerate}
As a consequence any germ $f$ of $T\setminus \Delta^{(r)}$ is topologically finite-determined with respect to $\tilde{I}$.
Furthermore, if $\KK=\CC$, given any open subset of $W'\subset W$, the filtration by succesive discriminats and the topological and finite-determinacy 
cummulative residual discriminants for $T_{W'}$ are the restrictions over $W'$ of the corresponding objets for $T$.}\\

This theorem shows, in particular, that, given any ideal sheaf of analytic functions, the subset of functions that are not topologically finite-determined 
with respect to it is very small (we can think of it as an infinite-codimension subset). In contrast with Theorem~\ref{old} we can not provide uniform
finite-determinacy bounds for prescribed codimension, that is, we can not ensure that for a prescribed integer $i$, there is another positive integer $r$ for 
which the subset $A_i$ necessarily $r$-determined. This is because 
in the jet-spaces $J^r(X,\tilde{I})$, instead of the algebraic structure present in ordinary jet-spaces, we have just an analytic structure. Furthermore, 
the subsets in which $\calO_{\CC^n,O}$ decomposes according to Varchenko's Theorem are invariant by the whole 
$\calD_O$; as we want to work in a neighbourhood of the origin, rather that just at the origin itself, we need to replace this by the property of beeing
flow-invariant. Nevertheless we have the following: 

\begin{remark}
\label{rema}
{\em
If we restrict to work at the origin, that is, to use the space $I=\tilde{I}_O$ instead of $J^{\infty}(X,\tilde{I})$, the corresponding jet spaces
$I/\mm^{r+1}_O\cap I$ are affine spaces. Then, if we assume $I$ to be generated by Nash functions, the arguments of this paper 
can be modified so that, if the starting subset $T$ of the Main Theorem is finitely determined and algebraic, then,
the subsets $A_i$ are finitely-determined and algebraic. Finally, let $\calD_{\tilde{I},O}$ 
be the subgroup of $\calD_O$ formed by diffeomorphisms which preserve the ideal $\tilde{I}_O$; if $T$ is assumed to be 
$\calD_{\tilde{I}_O}$-invariant, then the subsets $\{A_i\}_{i\geq 0}$, $\Gamma^{(r)}$ and $\Delta^{(r)}$  can be constructed
to be $\calD_{\tilde{I}_O}$-invariant.}
\end{remark}

\begin{proof}[Proof of the Main Theorem]
We show first the existence and uniqueness of the required objects satisfying all the requirements except Property~4. 
For any non-negative integer $j$ there is a unique filtration
\[T=A_0\supset ...\supset A_j\]
by closed $(C)$-analytic subsets, and two unique filtrations 
\[\Gamma_0\supset ...\supset\Gamma_j\]
\[\Delta_0\supset ...\supset\Delta_j\]
by closed subsets with the following properties:
\begin{itemize}
\item For any $i<j$ the sets $A_{i+1}\cup (\Gamma_j\cap A_i)$ and $A_{i+1}\cup (\Delta_j\cap A_i)$ are 
respectively the topological and finite-determinacy discriminants of $A_i$, beeing the set $A_i$ the analytic part in both cases.
\item The sets $\Gamma_j\cap A_j$ and $\Delta_j\cap A_j$ are respectively the topological and finite-determinacy discriminants of $A_j$.
\end{itemize}
The construction is obvious for $j=0$. Supposed that the filtrations have been constructed for a certain $j$, it is clear that $A_{j+1}$ must be defined to
be equal to the non-residual part of the topological discriminant of $A_j$, that is 
\[A_{j+1}:=(\Gamma_j\cap A_j)^{(a)}.\] 
The set $\Gamma_{j+1}$ only can be defined to be the union of the residual parts of the topological discriminants of $A_i$, for 
any $i\leq j$, with the whole topological discriminant of $A_{j+1}$. The definition of $\Delta_{j+1}$ is analogous.
We have shown by induction that the required filtrations can be constructed for any non-negative integer $j$ and are unique. It is easy to show that the 
infinite filtration 
\[A_0\supset ...\supset A_i\supset ...\]
and the closed sets $\Gamma^{(r)}:=\cap_{j\in\NN}\Gamma_j$ and $\Delta^{(r)}:=\cap_{j\in\NN}\Delta_j$ satisfy Properties $1-3$ from the statement of the 
Theorem. If $\KK=\CC$, using Lemma~\ref{restrabierto} it is easy to check that the filtration by succesive discriminats and the topological and 
finite-determinacy cummulative residual discriminants satisfy the compatibility condition concerning restrictions to open subsets of $W$.

It only remains to prove Property~4 when $T$ is flow-invariant. We only prove the statement concerning $A_i$, beeing the one concerning the $k$-determinacy 
discriminant analogous. We work by induction on $i$. Suppose that $A_j$ is flow-invariant for any $j\leq i$. We show that $A_{i+1}$ is flow-invariant. For this
we show that each irreducible component $C$ of $A_{i+1}$ is flow invariant.

Consider a increasing sequence $\{V_k\}_{k\in\NN}$ of open subsets of $U$ such that the closure $\overline{V}_k$ is compact and contained in $V_{k+1}$ and the
union $\cup_{k\in\NN}V_k$ equals $W$. Denote by $\partial V_k$ the boundary $\overline{V}_k\setminus V_k$. Define 
$d_k:=d(\overline{V}_{k},\partial V_{k+1})$, that is, the minimal euclidean distance between points of $\overline{V}_{k}$ and $\partial V_{k+1}$. 
Consider a decreasing sequence $\{\epsilon_k\}_{k\in\NN}$ of positive real numbers such that $\epsilon_k<d_k$ for any $k\in\NN$. 

Let $\calA_{k}$ the set of analytic diffeomorphisms $\phi$ defined on any neighbourhood $U$ of $\overline{V}_{k}$, satisfying $\phi(\overline{V}_k)\subset W$, 
and such that there exists a vector field $\theta\in\Gamma(W,\Theta_{\tilde{I},e})$ and a flow $\psi:U\times (a,b)\to W$ integrating $\theta$, with the 
property that $\phi=\psi_t$ for a certain $t\in (a,b)$. For any $\phi\in\calA_k$ we define $\delta(\phi):=\max\{||\phi(x)-x||:x\in\overline{V}_{k}\}$. 
Given any positive $\eta$ we define 
\[\calA_k^\eta:=\{\phi\in\calA_k:\delta(\phi)<\eta\}.\]
If $\eta\leq\epsilon_{k}$ then it is easy to show that $\phi(V_{k+1})\supset\overline{V}_{k}$ for any $\phi\in\calA_{k+1}^\eta$.
 
Let $r$ be the determinacy degree of the irreducible component $C$. Define
\[D_k(\eta):=\bigcap_{\phi\in\calA_{k+1}^\eta}\phi_*(C_{|V_{k+1}})_{|V_{k}},\]
which is a $r$-determined closed ($C$)-analytic subset of $J^\infty(X\cap V_{k},\tilde{I})$, for beeing intersection of such type of subsets.
Let $\{C_j\}_{j\in J}$ be the irreducible components of $A_{i+1}$ different from $C$; define 
\[A'_{i+1,k}(\eta):=D_k(\eta)\bigcup (\bigcup_{j\in J} C_{j|V_k}).\] 
Consider two germs $f$ and $g$ in the same path-connected component of 
\[A_{i|V_{k}}\setminus (A'_{i+1,k}(\eta)\cup \Gamma^{(r)}_{|V_k}).\]
We claim that $f$ and $g$
are topologically equivalent. Define 
\[X:=\cap_{\phi\in\calA_{k+1}^\eta}\phi_*(A_{i+1}\cup\Gamma^{(r)})_{|V_{k}}.\] 
As $X\subset A'_{i+1,k}(\eta)\cup\Gamma^{(r)}_{|V_k}$ the germs $f$ and $g$ are in the same path-connected component of $A_{i|V_{k}}\setminus X$. Let 
$\gamma:[0,1]\to A_{i|V_{k}}\setminus X$ be a continous path joining them. For each $t\in [0,1]$ there exists $\phi\in\calA_{k+1}^\eta$ such that 
$\gamma(t)$ is not in the closed subset $\phi_*(A_{i+1}\cup\Gamma^{(r)})_{|V_{k}}$; therefore, there exists a positive $\xi$ such that 
$\gamma(t-\xi,t+\xi)$ does not meet $\phi_*(A_{i+1}\cup\Gamma^{(r)})_{|V_{k}}$. Hence $\phi^*(\gamma(t-\xi,t+\xi))$ is
included in $A_{i}\setminus\ (A_{i+1}\cup\Gamma^{(r)})$ and consequently all the germs of $\gamma(t-\xi,t+\xi)$ are topologically equivalent. Covering 
$[0,1]$ by intervals like $(t-\xi,t+\xi)$ we conclude the proof of our claim.

Now we show that $D_k(\eta)$ is flow invariant. We do it in two steps: 

\textbf{Step 1}: we prove that $D_k(\eta)$ is flow-invariant with respect to vector fields defined on the whole $W$. Consider such a 
vector field $\theta\in\Gamma(W,\Theta_{I,e})$. Choose a point $x\in V_{k}$, consider a neighbourhood $U$ of $x$ in $W$ and a flow 
$\psi:U\times (a,b)\to V_{k}$ (with $0\in (a,b)$) obtained by integration of $\theta$. We have to show that $\psi_{t*}D_k(\eta)_x=D_k(\eta)_{\psi_t(x)}$ for 
any $t\in (a,b)$.

To get lighter notation denote by $D'(\eta)$ the projection $\pi^\infty_{r}(D_k(\eta))$. As $D_k(\eta)$ is $r$-determined it is enough to show that 
$j^{r}\psi_{t*}D'(\eta)_x\subset D'(\eta)_{\psi_t(x)}$ for any $t\in (a,b)$, or, which is the same, that $j^{r}\psi_{t*}(f)$ belongs to 
$D'(\eta)$ for any $f\in D'(\eta)_x$. Choose $f\in D'(\eta)_x$; as $D'(\eta)$ is closed, the subset $L_f\subset (a,b)$ consisting of $t$'s such that 
$j^{r}\psi_{t*}(f)$ belongs to $D'(\eta)$ is a non empty (as $0$ belongs to it) closed subset. If we prove that $L_f$ is also open then we conclude by 
connectedness of $(a,b)$.

Consider $t\in L_f$, we need to find a neighbourhood of $t$ included in $L_f$; as $\psi_{t_1}\comp\psi_{t_2}=\psi_{t_1+t_2}$, we can assume without loosing 
generality that $t=0$. By noetherianity of germs of analytic 
subsets there exists a neighbourhood $N_f$ of $f$ in $J^r(X,\tilde{I})$ and a finite subset $\calA (f)$ of $\calA_{k+1}^\eta$ such that 
\[D'(\eta)\cap N_f=[\bigcap_{\phi\in\calA (f)}j^{r}\phi_*(\pi^\infty_{r}(C_{|V_{k+1}}))_{|V_{k}}]\bigcap N_f;\]
By the finiteness of $\calA (f)$, the number $\nu:=\max\{\phi\in\calA (f):\delta(\phi)\}$ is strictly smaller than $\eta$. 
Consider the compact subset $K:=\cup_{\phi\in\calA (f)}\phi(\overline{V}_{k+1})\subset W$. There exists a positive $\xi$ such that if $|t|<\xi$, the domain of 
definition of $\psi_t$ contains an open neighbourhood of $K$ and 
\[\max\{||(\psi_t)(x)-x||:x\in K\}<\eta-\nu.\] 
Given any $\phi\in\calA (f)$ and $t$ with $|t|<\xi$, the domain of definition of the composite $\psi_t\comp\phi$ is clearly 
a neighbourhood of $\overline{V}_{k+1}$. Moreover 
\[||\psi_t\comp\phi(x)-x||\leq ||\psi_t(\phi(x))-\phi(x)||+||\phi(x)-x||<\eta-\nu+\eta=\eta,\]
and hence $\psi_t\comp\phi$ belongs to $\calA_{k+1}^\eta$. 
Then we have 
\[j^r\psi_{t*}(D'(\eta)\cap N_f)=j^r\psi_{t*}\big([\bigcap_{\phi\in\calA (f)}j^{r}\phi_*(\pi^\infty_{r}(C_{|V_{k+1}}))_{|V_k}]\bigcap N_f\big)\subset\]
\[\subset\bigcap_{\phi\in\calA (f)}j^{r}(\psi_t\comp\phi)_*(\pi^\infty_{r}(C_{V_{|k+1}}))_{|V_k}\subset\bigcap_{\phi\in\calA_{k+1}^\eta}j^{r}(\phi)_*(\pi^\infty_{r}(C_{V_{k+1}}))_{|V_k}=D'(\eta)\]
for any $t\in(-\xi,\xi)$; therefore $(-\xi,\xi)$ is included in $L_f$. 

\textbf{Step 2}: Let $\theta_1$,...,$\theta_l$ be vector fields generating the sheaf $\Theta_{\tilde{I},e}$ over $W$. Consider an open subset $U\subset W$ and 
a section $\theta\in\Gamma(U,\Theta_{\tilde{I},e})$. For any $x\in U$ there exists a neighbourhood $U_x$ of $x$ in $U$ and analytic functions $g_1,...,g_h$ 
such that $\theta=\sum_{i=1}^hg_i\theta_i$. Choose local coordinates $(z_1,...,z_n)$ around $x$ with the property that each $z_i$ is defined on the whole 
$\KK^n$; perhaps having to shrink $U_x$ we can assume 
that the power series expansions of $g_1,...,g_h$ with respect to these coordinates are convergent on the whole $U_x$. 
Denote by $g_i^{(l)}$ the trucation of the power series expansion of $g_i$ at the $l$-th term; the functions $g_i^{(l)}$ are polynomials in the local 
coordinates, and hence their domain of definition is also $\KK^n$. Therefore for each positive integer $l$ we can define the vector field 
$\theta^{(l)}\in\Gamma(W,\Theta_{\tilde{I},e})$ by the formula $\theta^{(l)}=\sum_{i=1}^hg_i^{(l)}\theta_i$. The sequence of vector fields 
$\{\theta^{(l)}\}_{l\in\NN}$ converges to $\theta$ on $U_x$.

The following statement can be deduced easily from the continous dependence of the solutions of differential equations with respect to a parameter (see 
\cite{AAA}~Ch.~1,~\S~2.8): there exists a positive integer $N$, an open neighbourhood $U'_x\subset U_x$ of $x$ containing the origin, and a positive real 
number $\xi$ such that 
\begin{enumerate}
\item For any $l\geq N$ there exists a flow $\psi^{(l)}:U'_x\times (-\xi,\xi)\to U_x$ integrating the vector field $\theta^{(l)}$.
\item There exists a flow $\psi:U'_x\times (-\xi,\xi)\to U_x$ integrating the vector field $\theta$, and sequence of mappings 
$\{\psi^{(l)}\}_{l\geq N}$ converges to $\psi$.  
\end{enumerate}

Suppose that we are given any flow $\psi:U'\times (a,b)\to U$ (with $0\in (a,b)$ integrating $\theta$. Define $L_f$ as in Step~1. Again it is enough to show 
that $L_f$ is open. Using that $\psi_{t_1+t_2}=\psi_{t_1}\comp\psi_{t_2}$ the proof can be reduced to the existence of $\xi>0$ such that $(-\xi,\xi)$ is 
contained in $L_f$. Choose a positive $\xi$ and an open neighbourhood $U'_x\subset U_x$ of $x$ such that the flow $\psi:U'_x\times (-\xi,\xi)\to U$ a limit of 
flows $\psi^{(l)}:U'_x\times (-\xi,\xi)\to U$ integrating vector fields $\{\theta^{(l)}\}_{l\in\NN}$ defined on the whole $W$. Then $j^r\psi_{t*}(f)$ is the 
limit of the sequence $\{j^r\psi^{(l)}_{t*}(f)\}_{l\in\NN}$. For beeing $\theta^{(l)}$ defined on the whole $W$, by Step~1, we have 
$j^r\psi^{(l)}_{t*}(f)\in D'(\eta)$ for any $l\in\NN$ and any $t\in(-\xi,\xi)$; then, as $D'(\eta)$ is closed, we have $j^r\psi_{t*}(f)\in D'(\eta)$ for any 
$t\in(-\xi,\xi)$ . This concludes the proof of
the flow-invariance of $D_k(\eta)$.

We claim that the restriction $D_{k+1}(\epsilon_{k+1})_{|V_{k}}$ is equal to $D_k(\epsilon_{k})$. As $V_{k+1}\subset V_{k+2}$ and 
$\epsilon_{k+1}<\epsilon_k$ we 
have $\calA_{k+2}^{\epsilon_{k+1}}\subset\calA_{k+1}^{\epsilon_k}$. Consequently $D_{k+1}(\epsilon_{k+1})_{|V_{k}}\supset D_k(\epsilon_{k})$. 
Obviously 
\[\bigcap_{\phi\in\calA_{k+1}^{\epsilon_{k}}}\phi_*(D_{k+1}(\epsilon_{k+1}))_{|V_{k}}\subset\bigcap_{\phi\in\calA_{k+1}^{\epsilon_{k}}}\phi_*(C_{|V_{k+1}})_{|V_{k}}=D_k(\epsilon_{k}).\]
On the other hand by the flow-invariance of $D_{k+1}(\epsilon_{k+1})_{|V_{k}}$  the first term of the last expression is equal to 
$D^{\epsilon_{k+1}}_{k+1|V_{k}}$. This shows our claim. It follows that the union
\[D:=\bigcup_{k\in\NN} D_k(\epsilon_{k})\]
is a $r$-determined closed $(C)$-analytic subset of $C$ which is flow-invariant. We define $A'_{i+1}$ as the union of $D$ with all the irreducible components
of $A_{i+1}$ different from $C$. Any two germs in the same connected component of $A_{i}\setminus (A'_{i+1}\cup\Gamma^{(r)})$ are topologically equivalent 
(it is enough to check this statement at the restriction over each $V_k$, and this has been already shown). Therefore the set 
$A'_{i+1}\cup(\Gamma^{(r)}\cap A_i)$ contains the topological discriminant of $A_i$. Taking analytic parts we get $A'_{i+1}\supset A_{i+1}$, which implies that
$D=C$. Consequently $C$ is flow-invariant.
\end{proof}

\section{Generalization of Varchenko's method} 

The overall structure of the proof of Proposition~\ref{Varchenko} follows~\cite{Va1}; it is based on an algorithm that shows the existence of a generic 
R-L-topological type for any family of functions, and a definition of 
so-called optimal germs in each finite-determined locally closed analytic subset of 
$J^\infty(X,\tilde{I})$. Nevertheless a straightforward generalization of Varchenko's proof to our case
does not work, mostly beacuse in his definition 
of optimal germs it is needed to perform certain modifications of functions that 
would take us outside the ideal sheaf we are working with, and also because he works with germs at the origin and we want to deal with a neighbourhood of it; 
this forces us to perform non-trivial modifications both in the algorithm and in the selection of optimal germs. 

For notational covenience we recall \S1.1~of~\cite{Va1}: let $\CC_n[z]$ be the space of monic polynomials 
$z^n+a_{n-1}z^{n-1}+...+a_0$; it is an affine space whose coordinates are $a_0,...,a_n$. For each sequence of 
positive integers $i_1,...,i_k$ such that $i_1+...+i_k=n$ we consider the subset of $\CC_n[x]$ consisting of polynomials 
with $k$ roots of multiplicities $i_1,...,i_k$; this defines a stratification of $\CC_n[z]$, whose strata will be called 
{\em multiplicity strata}. For each $m$ we consider $S_m$, the union of multiplicity strata containing polynomials
with less than $m$ different roots; the set $S_m$ is determined by a finite set of polynomial equations (with real coefficients) in $a_0,...,a_n$.  

Let $U\subset\KK^l$ be an open subset. An {\em $U$-family of functions} is, by definition, a $\KK$-valued analytic function $F$ defined on an open 
neighbourhood $V$ of $\{O\}\times U\subset\KK^n\times\KK^l$. Let $T\subset U$ be a closed ($C$)-analytic analytic subset, a {\em $T$-family of functions} 
is the retriction to $V\cap (\KK^n\times T)$ of a $U$-family of functions. With any $T$-family of  
functions of we associate its graph $\Gamma\in [V\cap (\KK^n\times T)]\times\KK$, i.e., the subvariety defined by the 
function $P^F:=u-F$, where $u$ is the coordinate function of the target $\KK$. 

Suppose $\KK=\RR$; view $\RR^n$ as the subset of points in $\CC^n$ with real coordinates. Let $F$ be a $T$-family of functions. As $T$ is $C$-analytic there 
exists an open subset $U^*\subset\CC^l$ and a closed analytic subset $T^*\subset U^*$ such that $U^*\cap\RR^n=U$ and $T^*\cap U=T$. According with 
Proposition~16~of~\cite{Na},~page~105, the subsets $T^*$ and $U^*$ can be chosen minimal in the following sense: if $T'$ is any other complex analytic subset 
of a neighbourhood of $T$ in $\CC^n$ for which $T'\supset T$, then $T'\cap W\supset T^*\cap W$ for a certain neighbourhood $W$ of $T$ in $\CC^n$. As $F$ can be
expressed 
locally as convergent power series, we can shrink $U^*$ so that there is a neighbourhood $V^*$ of $V$ in $\CC^n\times\CC^l$, containing $\{O\}\times U^*$, and 
a complex analytic function $F^*$ defined on $V^*$ whose restriction to $V$ is $F$. If $T^*$ is 
chosen minimal we say that the $T^*$-family defined by $F^*$ is a {\em minimal complex extension} of the $T$-family defined by $F$. We denote by $\Gamma^*$ 
the graph of
the $T^*$-family defined by $F^*$; clearly $\Gamma^*\cap (\RR^n\times\RR^l\times\RR)=\Gamma$.

Our aim is to show that generic functions of any $T$-family of functions are R-L-topologically equivalent. We work first for $\KK=\CC$, and then
explain the neccesary modifications needed for the real case.\\

\noindent
\textbf{Algorithm}: Fix a coordinate system $(x_1,...,x_n)$ of $\CC^n$. Consider a $T$-family of functions (for a certain $T\subset\CC^l$).
We describe an algorithm that constructs a new
coordinate system $(z_1,...,z_n)$ of $\CC^n$ (which will be called a {\em suitable coordinate system}), 
a proper analytic subset
$A\subset T$, positive continous functions $r_1,...,r_{n+1}$ defined over $T\setminus A$, non negative integers 
$k_1,k_2,...,k_{n+1}$ (with $k_1=0$) and a sequence of pseudopolynomials $P_1,...,P_{n+1}$ of the form $P_{n+1}:=u^{d_{n+1}}$, and
\[P_i(z_i,...,z_n,u,y)=z_i^{d_i}+\sum_{j=0}^{d_i-1}\alpha_i^j(z_{i+1},...,z_n,u,y)z_i^j,\]
for $i\leq n$, with $d_i\geq 0$ for any $i$, and $\alpha_i^j$ analytic in $U_{i+1}$, where 
\[U_i:=\{(z_i,...,z_n,u,y):y\in T\setminus A,|z_i|<r_i(y),...,|z_n|<r_n(y),|u|<r_{n+1}(y)\},\]
with the following properties: let $\Gamma_i:=V(P_i)\cup V(R_i)$, where $R_i:=u^{k_i}$, then
\begin{enumerate}
\item $\Gamma_1\cap U_1=\Gamma\cap U_1$.
\item For each $i\leq n$ the polynomials $P_i(z_i,a_{i+1},...,a_{n},b,c)$ are in the same multiplicity stratum as 
polynomials in $z_i$ if $(a_{i+1},...,a_{n},b,c)\in U_{i+1}\setminus\Gamma_{i+1}$.
\item The roots of the polynomial $P_i(z_i,a_{i+1},...,a_{n},b,c)$ are in the disc of radius $r_i(c)$ for any 
$(a_{i+1},...,a_{n},b,c)\in U_{i+1}$.
\item $\alpha_i^j(0,...,0,y)=0$ for any $i,j$.
\end{enumerate}
  
\begin{notation}
Let $f\in\CC\{z_1,...,z_n,u\}$, we define $\mult'(f):=\mult(f(z_1,...,z_n,0))$, and
$\wideg_{z_1}(f):=\mult(f(z_1,0,...,0))$.
\end{notation}

 Now we describe the algorithm under the assumption that $P^F(0,...,0,y)=0$ for any $y\in T$
(otherwise the existence of the claimed objects is easy):

\textbf{Step 1}: As $P^F(0,...,0,y)=0$ for any $y\in T$, and $u\!\not\!| P^F$, we 
deduce that $0<\mult'(P^F(\cdot,...,\cdot,y))<\infty$ for any $y\in T$; define
\[d_1:=\min\{\mult'(P^F(\cdot,...,\cdot,y)):y\in T\}.\]
Considering new coordinates $(z_1,z^1_2...,z^1_n)$ of $\CC^n$ related to $(x_1,...,x_n)$ by the formulas $x_1:=z_1$, 
$x_i:=z^i_1+\lambda^1_iz_1$ for $i>1$, where the $\lambda^1_i$'s are generic, we deduce that 
$\wideg_{z_1}(P^F(\centerdot,...,\centerdot,y))=d_1$ for certain $y\in T$. 

Define 
\[A_1:=\{y\in T:\wideg_{z_1}(P^F(\cdot,...,\cdot,y)>d_1\};\]
clearly $A_1$ is a proper analytic subset of $T$. 

By Weierstrass Preparation Theorem applied 
in a neighbourhood of the set 
\[\{(0,...,0,y):y\in T\setminus A_1\},\] 
there exists a neighbourhood $V_1$ of $T\setminus A_1$
in $\CC^l$ and positive continous functions $r_1,r^1_2,...,r^1_{n+1}$ defined on $V_1$ 
such that $P^F$ decomposes on the set 
\[\{(z_1,z^1_2,...,z^1_n,u,y):y\in V_1,|z_1|<r_1(y),...,|z^1_n|<r^1_n(y),|u|<r^1_{n+1}(y)\}\]
as $P^F=\phi P_1$ where $\phi$ is analytic and does not vanish at any point of the set and $P_1$ is a pseudopolynomial 
of the form 
\[P_1(z_1,z_2^1...,z^1_n,u,y)=z_1^{d_1}+\sum_{j=0}^{d_1-1}\alpha_1^j(z^1_{2},...,z^1_n,u,y)z_1^j,\]
such that its coefficients $\alpha_1^j$ are analytic on the set 
\[\{(z^1_2,...,z^1_n,u,y):y\in V_1,|z^1_2|<r^1_2(y),...,|z^1_n|<r^1_n(y),|u|<r^1_{n+1}(y)\}.\] 
Set $k_1=0$. Choosing $r^1_2,...,r^1_{n+1}$ small enough we ensure that Property~3 is satisfied for $P_1$.

\textbf{Step i} (for $1<i\leq n$): in the previous step we have constructed a system of coordinates 
$z_1,...z_{i-1},z_i^{i-1},...,z_n^{i-1}$ of $\CC^n$, a proper analytic
subset $A_{i-1}\subset T$, positive continous functions $r_1,...r_{i-1},r^{i-1}_i,r^{i-1}_{n+1}$, a neighbourhood $V_{i-1}$
of $T\setminus A_{i-1}$ in $\CC^l$ and a pseudopolynomial 
$P_{i-1}=z_{i-1}^{d_{i-1}}+\sum_{j=0}^{d_{i-1}-1}\alpha_{i-1}^jz_{i-1}^j$ such that the 
functions $\alpha_{i-1}^j$ are analytic on the set $U''_i$ defined by    
\[\{(z^{i-1}_i,...,z^{i-1}_n,u,y):y\in V_{i-1},|z_i^{i-1}|<r^{i-1}_i(y),...,|z^{i-1}_n|<r^{i-1}_n(y),|u|<r^{i-1}_{n+1}(y)\}.\]

Consider $P_{i-1}$ as a family of polynomials of $\CC_{d_{i-1}}[z_{i-1}]$ parametrized by the set $U'_{i}$ defined by
\[\{(z^{i-1}_i,...,y):y\in T\setminus A_{i-1},|z^{i-1}_i|<r^{i-1}_i(y),...,|z^{i-1}_n|<r^{i-1}_n(y),|u|<r^{i-1}_{n+1}(y)\}.\]
Since $T$ is irreducible there exists a multiplicity stratum of $\CC_{d_{i-1}}[z_{i-1}]$ whose closure contains all the 
polynomials of this family, and such that there is a polynomial of the family belonging precisely to this stratum.
Let this stratum contain polynomials with $m_i$ different roots; let $G_1,...,G_{k_i}$ be the polynomials in the variables
$a_0,...a_{d_{i-1}}$  determining the set $S_{m_i}$ in $\CC_{d_{i-1}}[z_{i-1}]$; we order them so that the first $s_i$ 
(for a certain positive $s_i$) of them are the polynomials that does not vanish identically in $U'_{i}$ when we substitute the 
variables $a_j$'s by the functions $\alpha_{i-1}^j$'s. Define an analytic function on the set $U''_i$ by the formula
\[P'_i:=\prod_{j=1}^{s_i}G_j(\alpha_{i-1}^0,...,\alpha_{i-1}^{d_{i-1}-1}).\]
Perhaps having to substitute $r^{i-1}_{n+1}$ by another smaller positive continous function, we can 
assume that $P'_i$ admits a unique expression in $U''_i$ as
\[P'_i=\sum_{k=0}^\infty\psi_ku^k,\]
where $\psi_k$ is an analytic function on $U''_i$ not depending on $u$. Let $k'_i$ be minimal such that the restriction of $\psi_{k'_i|U'_i}$ is not 
identically zero. Define $P''_i:=\sum_{k=k'_i}^\infty\psi_ku^k$; clearly $P''_{i|U'_i}=P'_{i|U'_i}$. Define $k_i:=k_{i-1}+k'_i$, and $P'''_i:=P''_i/u^{k'_i}$. 
Clearly we have 
\[d_i:=\min\{\mult'(P'''_i(·,...,·,y)):y\in T\setminus A_{i-1}\}<\infty.\]
We consider several cases:

\textsc{Case 1} ($d_i=0$): We choose the definitive coordinate system 
$(z_1,...,z_n)$ equal to $(z_1,...,z_{i-1},z^{i-1}_i,...,z_n^{i-1})$. Define the closed analytic proper subset $A$ as 
\[A:=A_{i-1}\cup\{y\in T\setminus A_{i-1}:\mult'(P'''_i(\cdot,...,\cdot,y)>0\};\] 
define 
$P_i=...=P_{n+1}=1$, take $k_{n+1}=...=k_{i+1}=0$ and choose, for any $j\geq i$ the positive continous function $r_j$ upper bounded by
$r^{i-1}_j$ and small enough so that the intersection $\{P''_i=0\}\cap U_i$ is empty. The algorithm concludes here.

\textsc{Case 2} ($d_i>0$): We consider new coordinates $(z_1,...,z_i,z_{i+1}^i,...,z_n^i)$ related to the previous ones by the formulas 
$z^{i-1}_i=z_i$ and $z^{i-1}_j=z^i_j+\lambda^i_jz_i$ for $j>i$. Choosing the $\lambda^i_j$'s generic we obtain that the set
\[A_i:=A_{i-1}\cup\{y\in T\setminus A_{i-1}:\wideg_{z_i}(P'''_i)>d_i\}\]
is a proper analytic subset of $T$. 

By Weierstrass Preparation Theorem applied 
in a neighbourhood of the set 
\[\{(0,...,0,y):y\in T\setminus A_i\},\] 
there exists a neighbourhood $V_i$ of $T\setminus A_i$
in $\CC^l$ and positive continous functions $r_i,r^{i}_{i+1}...,r^i_{n+1}$ such that $P''_i$ decomposes on the set 
\[\{(z_i,z^{i}_{i+1},...,z^i_n,u,y):y\in V_i,|z_i|<r_i(y),...,|z^i_n|<r^i_n(y),|u|<r^i_{n+1}(y)\}\]
as $P''_i=\phi P_i$ where $\phi$ is analytic and does not vanish at any point of the set, and $P_i$ is a pseudopolynomial of
the form
\[P_i(z_i,z^i_{i+1}...,z^i_n,u,y)=z_i^{d_i}+\sum_{j=0}^{d_i-1}\alpha_i^j(z^i_{i+1},...,z^i_n,u,y)z_i^j,\]
such that its coefficients $\alpha_i^j$ are analytic on 
\[\{(z^{i}_{i+1},...,z^i_n,u,y):y\in V_i,|z^i_{i+1}|<r^i_{i+1}(y),...,|z^i_n|<r^i_n(y),|u|<r^i_{n+1}(y)\}.\] 
Choosing $r^i_{i+1},...,r^i_{n+1}$ small enough we ensure that Property~3 is satisfied for $P_i$.

\textbf{Step n+1}: this step runs parallel to the induction step ({\em Step i}). As 
$u$ is the only variable of $P'_{i+1}$ we are forcedly in Case 1. This concludes the 
algorithm.\\

\textbf{The real case}: Suppose now that $\KK=\RR$. Consider $T\subset U$, a closed $C$-analytic subset contained in an open subset of $\RR^n$; let 
$V$ be an open neighbourhood of $\{O\}\times U$ in $\RR^n\times\RR^l$; let 
$F:V\to\RR$ define a $T$-family of functions. Consider $T^*\subset U^*\subset\CC^n$ and a complex analytic function $F^*:V^*\to\CC$ (where $V^*$ is a 
neighbourhood of $\{O\}\times U^*$ in $\CC^n\times\CC^l$) such that the $T^*$-family given by $F$ is a minimal complex extension of the $T$-family defined by 
$F$. It is easy to check that the preceding algorithm can be applied to the $T^*$-family $F^*$ so that:
\begin{enumerate}
\item The initial coordinate system $(x_1,...,x_n)$ is a real coordinate system when restricted to $\RR^n$. In each step we choose a change of coordinates with
{\em real} matrix, so that the new coordinates functions $(z_1,...,z_n)$ for $\CC^n$ form are real coordinate system when restricted to $\RR^n$.
\item The subset $A\cap\RR^n$ is a proper $C$-analytic subset of $T$.
\item For any $i\leq n$ and $j\leq d_i-1$ the function $\alpha_i^j$ assumes real values when $z_i,...,z_n,u,y$ are real.
\end{enumerate}

Given a $T$-family of functions $F$, for any $t\in T$ we will denote by $F_{|t}$ the germ at the origin given by the restriction of $F$ to $\KK^n\times\{t\}$.

\begin{prop}
\label{topeq}
Let $\KK$ be either $\RR$ or $\CC$. Let $F$ be a $T$-family of functions. Let $A\subset T$ be the subset constructed in the preceding algorithm.
For any $t,t'$ in the same connected component of $T\setminus A$, the functions $F_{|t}$ and $F_{|t'}$ are topologically R-L-equivalent.
\end{prop}
\begin{proof}
Suppose $\KK=\CC$; let $\Gamma$ be the graph of $F$; 
view $\Gamma$ as a family of analytic hypersurfaces parametrized by $T$. Proposition~3.1~of~\cite{Va1} can be 
adapted to apply (with minor changes in its proof) to this setting. If $\KK=\RR$ we consider a $T^*$-family $F^*$ which is a minimal complex extension of the 
$T$-family $F$; let $\Gamma^*$ be the graph of $F^*$. Then Proposition~3.2~of~\cite{Va1} can be adapted to apply (with changes in its proof) to 
the family of hypersurfaces $\Gamma^*$.    

Notice that we have designed our algorithm so that the matrix 
$(c_{i,j})$ relating the coordinate systems $(x_1,...,x_n,u)$ and $(z_1,...,z_n,u)$ of $\KK^n\times\KK$ has block form $c_{n+1,i}=c_{i,n+1}$ for any 
$i\neq n+1$. Since we are able to adapt Propositions~3.1~and~3.2~of~\cite{Va1}, the proofs of~Propositions~4.1~and~4.2~of~\cite{Va1} apply word by word in our 
case. Applying them respectively in the complex and real case, our result follows.
\end{proof}
 
Let $X$ be any $\tilde{I}$-stratum of $W$; let $T\subset J^\infty(X,\tilde{I})$ be an irreducible ($C$)-analytic subset. Before proving 
Proposition~\ref{Varchenko} we have to distinguish a special class of 
germs in $T$, which will be called {\em optimal germs}; this class is in a certain sense ``finitely determined and open''. 
We select the germs that we will call optimal in the following way:

\textbf{Search for Optimal Germs:} the differences between the Algorithm described above and Varchenko's algorithm forces us to introduce also some 
different features in the selection of optimal germs. In particular, Varchenko's search for optimal germs can be performed 
with {\em any coordinate system of $\KK^n$}; this will
not be the case in our situation: as the search for optimal germs advances we will need to modify our original coordinate
system, getting at the end a new one that will be regarded as {\em good coordinate system with respect to $T$}. We will 
proceed in several stages in which we will select smaller subsets of $T$ each time.

We fix an initial coordinate system $(x_1,...,x_n)$ for $\KK^n$; given any $f\in J^\infty(X,\tilde{I})$ we view it as a convergent power series of 
$\CC\{x_1,...,x_n\}$ by taking its Taylor expansion at $pr_\infty(f)$. 

\textbf{Stage 1}.
Given any $f\in T$ we consider $P^f=u-f\in\CC\{x_1,...,x_n,u\}$; obviously $\mult'(P^f)<\infty$, and 
therefore 
\[d_1:=\min\{\mult'(P^f):f\in T\}<\infty.\]
Considering $z_1,z_2^1,...,z_n^1$, a new coordinate system 
of $\KK^n$ related with the old one by formulas of the form $x_1=z_1$, $x_i=z^1_i+\lambda^1_iz_1$, with the 
$\lambda^1_i$'s real and generic enough, we obtain that $\min\{\wideg_{z_1}(P^f):f\in T\}=d_1$. Define the non-empty set
\[\calE^1_1:=\{f\in T/\wideg_{z_1}(P^f)=d_1\}.\] 

By Weierstrass Preparation Theorem , given any $f\in\calE^1_1$, it is possible to find positive numbers $r_1(f),...,r_{n+1}(f)$ such 
that $P^f$ can be decomposed over the open subset 
\[U^1_1(f):=\{(z_1,z^1_2,...,z^1_n,u):|z_1|\leq r_1(f),...,|u|\leq r_{n+1}(f)\}\]
as $P^f=\phi_1P_1[f]$, where $\phi_1$ does not vanish anywhere in $U_1$ and $P_1[f]$ is a pseudopolynomial of the form
\[P_1[f](z_1,...,u):=z_1^{d_1}+\sum_{j=1}^{d_1-1}z_1^j\alpha_1^j[f](z_2^1,...,u).\]
View $P_1[f]$ as a family of polynomials of $\CC_{d_1}[z_1]$ (even when $\KK=\RR$) parametrized by the open subset    
\[U^1_2(f):=\{(z^1_2,...,z^1_n,u):|z^1_2|\leq r_2(f),...,|u|\leq r_{n+1}(f)\}.\]
There exists a stratum $S_1(f)$ of $\CC_{d_1}[z_1]$ whose closure contains the whole family, and such that there is a 
member of the family belonging to it; let $m_1(f)$ be the number of roots of a generic element of $S_1(f)$. Define
$m_1^0:=\max\{m_1(f):f\in\calE^1_1\}$, and let $S_1^0$ be a stratum of  $\CC_{d_1}[z_1]$ such that there is 
$f\in\calE_1^1$ with $S_1(f)=S_1^0$ and $m_1(f)=m_1^0$; define the non-empty set
\[\calE^2_1:=\{f\in\calE_1^1/S_1(f):=S_1^0\}.\]

Let $z_1^{d_1}+\sum_{j=1}^{d_1-1}a_jz_1^j$ be a generic polynomial $\CC_{d_1}[z_1]$; there are polynomials $Q_1,...,Q_s$ (with real coefficients)
in the variables $a_0,...,a_{d_1-1}$ whose set of common zeros determines the set of polynomials in $\CC_{d_1}[z_1]$ with 
less than $m_1^0$ roots. Let $f\in\calE_1^2$, denote by $T_i[f]$ the analytic functions in $z_2^1,...,z_n^1,u$ obtained 
substituting in $Q_1$ the $a_j$'s by the $\alpha_1^j[f]$'s, define $s_1(f)$ to be the number of $T_j[f]$'s that do not 
vanish identically in $U_2(f)$; clearly $s_1(f)>0$ for any $f\in\calE_1^2$. Define $s_1^0:=\max\{s_1(f):f\in\calE_1^2\}$, 
choose $f\in\calE_1^2$ such that $s_1(f)=s_1^0$, up to a re-ordering we can assume that 
\[P'_2[f]:=\prod_{i=1}^{s_1^0}T_i[f]\] 
does not vanish identically in $U_2(f)$. Define the non-empty set
\[\calE^3_1:=\{f\in\calE_2^1/P'_2[f]_{|U_2(f)}\not\equiv 0\}.\] 

For any $f\in\calE_1^3$ we let $k_1(f)$ be the maximal power of $u$ which divides $P'_2[f]$; we define
$k_1^0:=\min\{k_1(f):f\in\calE_1^3\}$, 
\[\calE_1:=\{f\in\calE_3^1/k_1(f)=k_1^0\},\]
and, for any $f\in\calE_1$,
\[P''_2[f]:=P'_2[f]/u^{k_1^0}.\]

\textbf{Stage i} (for $1<i\leq n$): in the previous stage we have choosen a coordinate system 
$z_1,...,z_{i-1},z_i^{i-1},...,z_n^{i-1}$ of $\KK^n$ and set of germs $\calE_{i-1}$; moreover for each $f\in\calE_{i-1}$ 
we have given positive numbers $r_1(f),...,r_{n+1}(f)$ and an analytic function $P''_{i}[f]$ defined on 
\[U^{i-1}_{i}(f)=\{(z^{i-1}_i,...,z^{i-1}_n,u):|z^{i-1}_i|<r_i(f),...,|u|<r_{n+1}(f)\}\]
such that $u\not|P''_i[f]$ (hence $\mult'(P''_i[f])<\infty$). Consider 
\[d_i:=\min\{\mult'(P''_{i}[f]):f\in\calE_{i-1}\}.\]
Define $\{z_i,z^i_{i+1},...,z_n^i\}$ by the formulas $z^{i-1}_i=z_i$ and
$z^{i-1}_j=z^i_j+\lambda^i_jz_i$, with $\lambda_j^i$ real for $j>i$; for any $f\in\calE_{i-1}$ we express $P''_{i}[f]$ 
respect to the new variables $z_1,...,z_i,z^i_{i+1},...,z_n^i,u$. Then, if the $\lambda_j^i$'s are chosen generic enough, there 
exists $f\in\calE_{i-1}$ such that $\wideg_{z_i}(P''_{i}[f])=d_i$. Define the non-empty set
\[\calE_i^1:=\{f\in\calE_{i-1}:\wideg_{z_i}(P''_{i}[f])=d_i\}.\]

For each $f\in\calE_i^1$ we can diminish the numbers $r_i(f),...,r_{n+1}(f)$ so that the function
$P''_{i}[f]$ can be decomposed over the open subset 
\[U^i_1(f):=\{(z_i,z^i_{i+1},...,z^i_n,u):|z_i|\leq r_i(f),...,|u|\leq r_{n+1}(f)\}\]
as $P''_{i}[f]=\phi_iP_i[f]$, where $\phi_i$ does not vanish anywhere in $U^i_i$ and $P_i[f]$ is a pseudopolynomial of the 
form
\[P_i[f](z_i,...,u):=z_i^{d_i}+\sum_{j=1}^{d_i-1}z_i^j\alpha_i^j[f](z_{i+1}^i,...,u).\]
View $P_i[f]$ as a family of polynomials of $\CC_{d_i}[z_i]$ parametrized by the open subset    
\[U^i_{2}(f):=\{(z^i_{i+1},...,z^i_{i+1},u):|z^1_{i+1}|\leq r_{i+1}(f),...,|u|\leq r_{n+1}(f)\};\]

By analogy with Stage~1 we define numbers $m_i^0$, $s_i^0$, $k_i^0$, a  stratum $S_i^0\in\CC_{d_i}[z_i]$, and a decreasing
sequence of subsets $\calE_i^1\supset\calE_i^2\supset\calE_i^3\supset\calE_i$; moreover for each function $f\in\calE_i$ 
we construct functions $P'_{i+1}[f]$ and $P''_{i}[f]$ analytic in $z^i_{i+1},...,u$, such that $P'_{i+1}[f]=u^{k_i^0}P''_{i}[f]$ and 
$u\not|P''_{i}[f]$.

\textsc{Stage $n+1$}: for any $f\in\calE_n$ the function $P''_n[f]$ is a unit in $\KK\{u\}$. Moreover in the previous stage
we have constructed a coordinate system $(z_1,...,z_n)$ of $\KK^n$ which will be said to be {\em a good coordinate system 
with respect to} $T$. We will define the subset 
\begin{equation}
\label{optimo}
\calE_T\subset\pi_n^{-1}(T)
\end{equation}
 of {\em optimal germs of} $\pi_n^{-1}(T)$ 
{\em with respect to the coordinate system} $(z_1,...,z_n)$ as $\calE_T:=\calE_n$. This finishes the Search for Optimal 
Germs.

\begin{proof}[Proof of Proposition~\ref{Varchenko}]
We will start by proving that beeing optimal respect to a fixed good coordinate system $(z_1,...,z_n)$ is a finitely determined 
property, i.e., that there exists $s>r$ with the following property: given any $x\in X$ and any $g\in\tilde{I}_x\cap\mm^{s+1}_x$, a germ $f\in T_x$ belongs to
$\calE_T$ if and only if $f+g$ belongs to $\calE_T$. By convenience of the reader we repeat the statement of Proposition~4.2~of~\cite{Va1}:\\

\noindent
($\dagger$){\em For any two natural numbers $k$ and $p$ there exists third one $L(k,p)$ such that for any $f\in\CC\{z_1,...,z_n\}$ 
with $\wideg_{z_1}(f)=k$ and $g\in\mm^{L(k,p)}$ the following property holds: consider analytic functions $\phi$, $\phi'$, 
$P$, $P'$, such that $\phi(O)\neq 0$, $\phi'(O)\neq 0$ and $f=P\phi$, $f+g=P' \phi'$ in some neighbourhood of $O$, where 
$P$ and $P'$ are Weierstrass polynomials in $z_1$ of degree $k$; then $\phi-\phi'$ and $P-P'$ are in $\mm^p$.}\\

In order to choose $s$ we define the following numbers recursively:  
\[a_n:=k_{n-1}^0+\max\{d_n+1,L(d_n,d_n+k_n^0)\},\]
\[a_i:=k^0_{i-1}+\max\{d_i+1,L(d_i,a_{i+1}+d_i)\}\quad\text{for}\quad 2\leq i\leq n-1,\]
\[a_1:=\max\{d_1+1,L(d_1,a_2+d_1)\}.\]
Fix $s:=a_1$. With this choice, taking into account ($\dagger$) along the procedure of Search for Optimal Germs, it is 
easy to show that given $f\in T_x$ and $g\in\tilde{I}_x\cap\mm^{s+1}_x$ then $f\in\calE_T$ if and only if $f+g\in\calE_T$. 

Now we suppose that $f\in\calE_T$; we will show that then $f$ and $f+g$ are topologically equivalent. Consider a coordinate function $y$ for the affine line 
$\AAA^1_\KK$, and define the $\AAA^1_\KK$-family $F:=f+yg$. 
The germ $F_{|y}$ is optimal respect to $(z_1,...,z_n)$ for any $y\in\AAA^1_\KK$ for beeing $yg\in\tilde{I}_x\cap\mm^{s+1}_x$. 
It is easy to check that the Algorithm can be applied to 
the family $F$ taking as initial coordinate system $(z_1,...,z_n)$ and having in each step the {\em trivial} coordinate 
change; denote by $A$ the subset of $\AAA^1_\KK$ contructed in the algorithm. The Search for Optimal Germs has been designed in a compatible way with the 
Algorithm so that $F_{|y}$ is optimal with respect to $(z_1,...,z_n)$ if and only if $y\in\AAA^1_\KK\setminus A$. Therefore $A=\emptyset$. Applying 
Proposition~\ref{topeq} we obtain the topological equivalence of $f=F_{|0}$ and $f+g=F_{|1}$.

We are ready to prove the statement of the proposition. Recall that $\tilde{I}$ is generated over $W$ by the functions $f_1,...,f_m$; given any $x\in W$ we 
denote by $f_{i,x}$ the Taylor expansion of $f_i$ at $x$. Consider the $\KK$-vector space $E$ of 
polynomials of degree bounded by $S$, let $\{g_1,...,g_N\}$ be a basis of $V$; the product $W\times V^m$ is an open subset of $\KK^{n+mN}$, and a point of it 
is represented by a $mN+1$-uple $(x,\lambda^1_1,...,\lambda^1_N,...,\lambda^m_1,...,\lambda^N_m)$, with $x\in W$ and $\lambda^i_j\in\KK$. 
Consider the $W\times V^m$-family defined by the unique analytic function $F$ in a neighbourhood of $\{O\}\times W\times V^m$ in $\KK^n\times W\times V^m$ 
satisfying  
\[F_{|(x,h^i_j)}=\sum_{i=1}^m\sum_{j=1}^Ng_jf_{i,x}\]
for any $(x,h^i_j)\in W\times V^m$. Associated with $F$ we have a natural analytic mapping
\[\psi:W\times V^m\to J^s(W,\KK)\]
wich assigns to $(x,h^i_j)$ the $s$-jet of the germ $F_{|(x,h^i_j)}$ viewed as a germ at $x$. 
As $T$ is a $r$-determined ($C$)-analytic subset of $J^\infty(X,\tilde{I})$ and $s>r$, the subset $T':=\psi_s^{-1}(\pi^\infty_s(T))$ is a closed ($C$)-analytic
subset of $(W\setminus\partial X)\times V^m$; therefore we can consider the $T'$-family of functions obtained by restriction of $F$. 

It is easy to check that $\psi(T')=\pi^\infty_s(T)$; therefore there exists $t'\in T'$ such that $\psi(t')$ is the $s$-jet of an optimal germ of $T$ with 
respect to the coordinate system $(z_1,...,z_n)$. The Algorithm and the Search for Optimal Germs have been designed so that if we take into account that
beeing optimal with respect to the fixed coordinate system is an $s$-determined property, then 
\begin{enumerate}
\item The Algorithm can be applied 
to the $T'$-family of functions $F$ choosing at each stage a the trivial coordinate change.
\item A point $t'\in T'$ belongs to the ($C$)-analytic proper subset
$A'\subset T'$ constructed in the Algorithm if and only if $\psi(t')$ is not the $s$-jet of an optimal germ of $T$ with respect to the fixed coordinate system.
\end{enumerate}
Therefore $A'$ is of the form $\psi^{-1}(A'')$ with $A''\subset J^s(X,\tilde{I})$. The fact that $A''$ is ($C$)-analytic is deduced easily from 
the facts that $\psi_{|X\times V^m}:X\times V^m\to J^s(X,\tilde{I})$ is an epimorphism of trivial analytic vector bundles, and that $A'$ is ($C$)-analytic.
Define $A$ as the $s$ determined ($C$)-analytic subset $A:=(\pi^\infty_s)^{-1}(A'')$. 

Consider $f,f'$ in the same connected component of $T\setminus A$. Then there exist $t,t'$ in the same connected component of $T'\setminus A'$ such that 
$\psi(t)=\pi^\infty_s(f)$ and $\psi(t')=\pi^\infty_s(f')$. By Proposition~\ref{topeq} we have that $F_{|t}$ and $F_{|t'}$ are topologically equivalent. 
For beeing $t,t'\not\in A$ the germs $F_{|t}$ and $F_{|t'}$ are optimal with respect with the coordinate system $(z_1,...,z_n)$. As $f$ and $f'$ have 
respectively the same $s$-jet that $F_{|t}$ and $F_{|t'}$, and optimal germs are topologically $s$-determined, we conclude the topological equivalence of $f$ 
and $f'$. 
\end{proof}

\end{document}